\documentclass[11pt]{article}
\usepackage{amsmath,amssymb,amsthm}
\usepackage{algorithm}
\usepackage{algorithmic}
\usepackage{authblk}
\usepackage{color}
\usepackage{hyperref} 
\usepackage[left=1in,top=1in,right=1in,bottom=1in,letterpaper]{geometry}
\usepackage{enumerate}   
\usepackage{graphicx}

\numberwithin{equation}{section}

\newtheorem{theorem}{Theorem}[section]

\newtheorem{definition}[theorem]{Definition}
\newtheorem{example}{Example}[section]

\newtheorem{remark}[theorem]{Remark}

\newcommand{\X}{\mathcal{X}}

\newcommand{\argmin}{\mathop{\rm argmin}}
\newcommand{\argmax}{\mathop{\rm argmax}}
\newcommand{\br}{\mathbb{R}}

\newcommand{\half}{\frac{1}{2}}

\newcommand{\st}{\mathrm{s.t. }}

\newcommand{\Fix}{\mathrm{Fix}}

\newcommand{\be}{\begin{equation}}
\newcommand{\ee}{\end{equation}}
\newcommand{\ba}{\begin{array}}
\newcommand{\ea}{\end{array}}
\newcommand{\bad}{\begin{aligned}}
\newcommand{\ead}{\end{aligned}}

\newcommand{\setword}[2]{\phantomsection #1\def\@currentlabel{\unexpanded{#1}}\label{#2}}

\newcommand{\LCal}{\mathcal{L}}

\newcommand{\prox}{\mathrm{prox}}

\newcommand{\bone}{\mathbb{I}}
\newcommand{\onebf}{\mathbf{1}}

\newcommand{\intdom}{\mathrm{int\ dom}}
\newcommand{\dom}{\textrm{dom } }
\newcommand{\bpm}{\begin{pmatrix}}
\newcommand{\epm}{\end{pmatrix}}

\begin{document}
\title{Bregman Douglas-Rachford Splitting Method\thanks{The three authors started to work on this problem back in 2020, and this draft was ready in Nov 2021. We spent the last four years trying to prove the convergence of the algorithms under the most general setting, but did not succeed. Despite of that, we believe that the current contributions are significant. So we decided to share the results with the community now. Since this draft is from 2021, we apologize that the references are not up to date. We will include more recent works in a later version of the paper.}}
\author{Shiqian Ma$^1$, Lin Xiao$^2$, Renbo Zhao$^3$} 
\affil{$^1$Rice University, $^2$, $^3$University of Iowa}
\date{First version: Nov 5, 2021 \\ This version: Sep 4, 2025\thanks{The results in Sections 1-5 were written by Shiqian Ma in 2020, and the results in the Appendix were written by Shiqian Ma in Nov 2021. These results formed the first version of the paper which was ready on Nov 5, 2021. Shiqian Ma polished the writing in Sep 2025, which gives the current version. }}
\maketitle

\begin{abstract}
    In this paper, we propose the Bregman Douglas-Rachford splitting (BDRS) method and its variant Bregman Peaceman-Rachford splitting method for solving maximal monotone inclusion problem. We show that BDRS is equivalent to a Bregman alternating direction method of multipliers (ADMM) when applied to the dual of the problem. A special case of the Bregman ADMM is an alternating direction version of the exponential multiplier method. To the best of our knowledge, algorithms proposed in this paper are new to the literature. We also discuss how to use our algorithms to solve the discrete optimal transport (OT) problem. We prove the convergence of the algorithms under certain assumptions, though we point out that one assumption does not apply to the OT problem. 
\end{abstract}

\section{Introduction}

This paper studies the celebrated Douglas-Rachford splitting method (DRS) that has a long history dating back to 1956 for solving variational problems arising from numerical PDEs \cite{Douglas-Rachford-56}. The DRS later on became a very popular method for finding zeros of the sum of maximal monotone operators: 
\be\label{sum-2-operator}
\mbox{ Find } x, \ \st, \ 0 \in  A(x) + B(x),
\ee
where $A: \br^n \to \br^n$ and $B: \br^n \to \br^n$ are two maximal monotone operators. The problem \eqref{sum-2-operator} is also known as the monotone inclusion problem. The solution set of \eqref{sum-2-operator} is denoted as $(A+B)^{-1}(0)$.  
The DRS has been studied by many researchers under various settings \cite{Lions-Mercier-79,Gabay-Mercier-1976,Gabay-83,Glowinski-LeTallec-89,Fortin-Glowinski-1983,Eckstein-thesis-89,Eckstein-Bertsekas-1992,Bauschke-Moursi-DR-2017,Borwein-DR-nonconvex-2011,Combettes-Pesquet-DR-2007,Svaiter-DR-weak-convergence-2011,Li-Pong-DR-nonconvex,He-Yuan-DRS-rate-2015,Ryu-DRS-unique,Ryu-DR-pathological,Liang-Fadili-2017,Ryu-Yin-book,Bui-Combettes-DR-converge-weakly,Censor-r-sets-DR-2019}. When $A$ and $B$ are normal cone operators, the problem \eqref{sum-2-operator} reduces to a feasibility problem which seeks for a point in the intersection of two sets. We refer to the recent survey paper \cite{lindstrom-survey-2018} for more details on DRS for feasibility problems. The DRS for solving \eqref{sum-2-operator} can be written as:
\begin{subequations}\label{sum-2-operator-DRS}
\begin{align}
x^{k}     & := J_{\gamma_k B}(z^k) \label{sum-2-operator-DRS-1}\\ 
y^k        & := J_{\gamma_k A}(2x^k-z^k) \label{sum-2-operator-DRS-2}\\ 
z^{k+1} & := z^k - x^k + y^k, \label{sum-2-operator-DRS-3}
\end{align}
\end{subequations}
where $J_{T} := (I+T)^{-1}$ is called the resolvent of the operator $T$, and $\gamma_k > 0$ is a parameter. 

A well-known application of the DRS is the so-called alternating direction method of multipliers (ADMM), which is usually applied to solving the following convex minimization problem: 
\be\label{admm-sum-2}
\min_{u,v} \ f(u) + g(v), \ \st, \ Mu + Nv = b, \ u\in\br^p, v\in\br^q,
\ee
where $b\in\br^m$, $M\in\br^{m\times p}$ and $N\in\br^{m\times q}$, and $f$ and $g$ are both proper, closed and convex functions. The ADMM for solving \eqref{admm-sum-2} updates the iterates as follows: 
\begin{subequations}\label{admm-sum-2-alg}
\begin{align}
u^{k+1} & := \argmin_u \ \LCal_\beta(u,v^k;w^k) \label{admm-sum-2-alg-1}\\ 
v^{k+1} & := \argmin_v \ \LCal_\beta(u^{k+1},v;w^k) \label{admm-sum-2-alg-2} \\
w^{k+1} & := w^k + \beta(Mu^{k+1} + Nv^{k+1} - b). \label{admm-sum-2-alg-3}
\end{align}
\end{subequations}
Here 
\be\label{admm-sum-2-AL-Fn}
\LCal_\beta(u,v;w) := f(u) + g(v) + \langle w, Mu + Nv - b\rangle + \frac{\beta}{2}\|Mu + Nv - b\|_2^2,
\ee
is the augmented Lagrangian function for \eqref{admm-sum-2}, where $w$ denotes the Lagrange multiplier, and $\beta> 0$ is a penalty parameter. As proved by Gabay \cite{Gabay-83}, ADMM \eqref{admm-sum-2-alg} for solving \eqref{admm-sum-2} is a special case of DRS \eqref{sum-2-operator-DRS} applied to solving the dual problem of \eqref{admm-sum-2}, whose optimality condition is in the form of \eqref{sum-2-operator}. More specifically, the dual problem of \eqref{admm-sum-2} is given by
\be\label{admm-sum-2-dual}
\min_x \ f^*(-M^\top x) + g^*(-N^\top x) + b^\top x,
\ee
where $f^*$ and $g^*$ are the conjugate functions of $f$ and $g$, respectively. The optimality condition of \eqref{admm-sum-2-dual} is: 
\be\label{admm-sum-2-dual-OptCond-0}
0\in -M\partial f^*(-M^\top x) - N \partial g^*(-N^\top x) + b,
\ee
which is in the form of \eqref{sum-2-operator} by defining $A(x) = -M\partial f^*(-M^\top x)$ and $B(x) = -N \partial g^*(-N^\top x) + b$. Appying DRS \eqref{sum-2-operator-DRS} to \eqref{admm-sum-2-dual-OptCond-0} gives the ADMM \eqref{admm-sum-2-alg}.

ADMM has received significant attention due to its applications in signal processing, image processing, semidefinite programming and statistics \cite{Combettes-Pesquet-DR-2007,Goldstein-Osher-08,Wang-Yang-Yin-Zhang-2008,Wen-Goldfarb-Yin-2009,Ma-Aybat-RPCA-survey}. It is not possible to exhaust the vast literature on ADMM, and we thus refer the reader to the following survey papers for more details on the theory and applications of ADMM and its variants \cite{Boyd-etal-ADM-survey-2011,Eckstein-Yao-tutorial-admm,Glowinski-admm-tutorial,Ma-Hong-admm-survey}. 

The efficiency of the DRS \eqref{sum-2-operator-DRS} relies on the assumption that both $J_{\gamma_k A}(z)$ and $J_{\gamma_k B}(z)$ can be computed easily, and similarly, the efficiency of the ADMM \eqref{admm-sum-2-alg} replies on the assumption that the two minimization subproblems are easy to solve. As we will discuss later, in certain applications, these computations are not easy (even when $M=N=I$) and more general Bregman distances need to be considered when designing these algorithms. However, to the best of our knowledge, Bregman DRS (BDRS) was not considered in the literature before. There exists one work on Bregman ADMM \cite{wang-bregman-admm}, but this algorithm only applies to some special class of problems and is different from the algorithms that we consider in this paper. 

{\bf Our contributions.} In this paper, we target to design the BDRS algorithm, and our main contributions are as follows. 
\begin{enumerate}[(i)]
\item We design the first BDRS algorithm in the literature, and analyze its connections with several existing methods. We also propose a Bregman Peaceman-Rachford splitting (BPRS) method, which is a close variant of BDRS. We show that BDRS is equivalent to a Bregman ADMM algorithm when applied to its dual.
\item We show that if the Bregman distance is generated by the Boltzmann-Shannon entropy, then when applied to the dual problem of linear inequality constrained convex programming problem, our BDRS gives an alternating direction version of the exponential multiplier method. We name this algorithm ADEMM, and this is also a new algorithm, to the best of our knowledge. 
\item We discuss how to use our BDRS and ADEMM to solve the discrete optimal transport (OT) problem, and discuss how they relate to and why they are better than the Sinkhorn's algorithm. 
\item We prove the convergence of BDRS under certain assumptions, though we want to point out that one of the assumptions does not apply to the OT problem. 
\end{enumerate}

{\bf Organization.} The rest of this paper is organized as follows. In Section \ref{sec:bregman-algs} we provide some preliminaries on Bregman distance and algorithms based on it. In Section \ref{sec:BDRS}, we propose our BDRS and ADEMM algorithms and discuss their connections with existing algorithms. In Section \ref{sec:BPRS}, we propose the BPRS algorithm. In Section \ref{sec:OT}, we discuss how to apply our proposed algorithms to solving the discrete optimal transport problem. We draw some conclusions in Section \ref{sec:conclusion}. In the appendix, we provide the convergence analysis for BPRS and BDRS, both under certain assumptions. 

\section{Preliminaries and Existing Bregman Algorithms}\label{sec:bregman-algs}
In this section, we briefly review the basics of Bregman distance and existing algorithms that use Bregman distance. The Bregman distance was first proposed by Bregman \cite{Bregman-67} in a primal-dual method for solving linearly constrained convex programming problems that involves non-orthogonal projections onto hyperplanes. This method was further studied by Censor and Lent \cite{Censor-Lent-1981}, and De Pierro and Iusem \cite{DePierro-Iusem-1986}. 

We now introduce the notions of Legendre function and Bregman distance. 
\begin{definition}[\cite{Rockafellar-book-70}]\label{def:Legendre-function}
A function $h$ is called a Legendre function, if it is proper, lower semicontinuous, strictly convex and essentially smooth.
\end{definition}
A Legendre function enjoys the following two useful properties:
\begin{enumerate}[(i)]
\item $h$ is Legendre if and only if its conjugate $h^*$ is Legendre.
\item The gradient of a Legendre function $h$ is a bijection from $\intdom h$ to $\intdom h^*$, and its inverse is the gradient of the conjugate, that is, we have $(\nabla h)^{-1} = \nabla h^*$.
\end{enumerate}
\begin{definition}\label{def:Bregman-distance}
For a Legendre function $h$, the Bregman distance corresponding to $h$ is defined as 
\be\label{def:Bregman-dist}
D_h(x,y) := h(x) - h(y) - \langle \nabla h(y), x-y\rangle.
\ee 
\end{definition}
The following Bregman distances are commonly seen in practice (for more examples, see \cite{Kiwiel-97}).
\begin{example}\label{Example:Bregman-distance}
\begin{enumerate}[(i)]
\item Energy: If $h(x) = \half\|x\|_2$, then $D_h(x,y) = \half \|x-y\|_2^2$ is the Euclidean distance. 
\item Quadratic: If $h(x) = \half x^\top Lx$ with matrix $L \succ 0$, then $D_h(x,y) = \half \|x-y\|_L^2 = \half (x-y)^\top L (x-y)$.
\item Boltzmann-Shannon entropy: If $h(x)$ is the Boltzmann-Shannon entropy function defined as $h(x) =  \sum_i x_i (\log x_i-1)$, then $D_h(x,y) = \sum_i x_i\log \frac{x_i}{y_i}-x_i+y_i$ is the Kullback-Leibler (KL) divergence. Note that the domain of $h$ is $\dom h = \br_{++}^n = \{x\mid x>0\}$. Moreover, $h^*(x) = \sum_i e^{x_i}$.
\item Burg's entropy: If $h(x) = -\sum_i \log x_i$, then $D_h(x,y) = -\sum_i \log\frac{x_i}{y_i}+\frac{x_i}{y_i}-1$. Note that the domain of $h$ is $\dom h = \br_{++}^n$.
\end{enumerate}
\end{example}

We now define a few useful Bregman operators. 

\begin{definition}[Bregman forward operator, Bregman resolvent operator, Bregman reflection operator, Bregman Mann's operator]\label{def-Bregman-operators}
Use $h$ to denote a Legendre function.
The Bregman forward operator for a single-valued operator $T$ is defined as 
\be\label{def-Bregman-forward-operator}
F_T^h := \nabla h^*\circ (\nabla h - T).
\ee
The Bregman resolvent operator of a maximal monotone operator $T$ is defined as \cite{Eckstein-93}
\be\label{def-Bregman-resolvent}
J_T^h := (\nabla h + T)^{-1}\circ \nabla h.
\ee
The Bregman reflection operator of a maximal monotone operator $T$ is defined as
\be\label{def-reflection-op}
R_T^h := \nabla h^*\circ (2\nabla h\circ J_T^h - \nabla h),
\ee
and the Bregman Mann's operator of a maximal monotone operator $T$ is defined as
\be\label{def-Mann}
M_\alpha^h(T) := \nabla h^*\circ(\alpha \nabla h + (1-\alpha)\nabla h\circ T),
\ee
\end{definition}
The notion of Bregman resolvent operator \eqref{def-Bregman-resolvent} was first proposed by Eckstein in \cite{Eckstein-93}. When $h(\cdot) = \half \|\cdot\|_2^2$, $J_T^h$ reduces to $J_T := (I+T)^{-1}$, because $\nabla h = I$. The definitions of Bregman relection operator \eqref{def-reflection-op} and Bregman Mann's operator \eqref{def-Mann} are new to the literature to the best of our knowledge.

We now review several classes of Bregman algorithms for both monotone inclusion and convex optimization problems. 

\subsection{Bregman Gradient Method and Mirror Descent Method} 

For unconstrained convex minimization problem 
\be\label{prob-uncon}\min_{x\in \X} \ f(x)\ee
where $f$ is convex and smooth and $\X \subseteq \br^n$ is a convex set, Nemirovski and Yudin \cite{ny83} proposed the mirror descent algorithm which iterates as
\be\label{BGM-alg-explicit}
x^{k+1} := \nabla h^*(\nabla h(x^k) - \gamma_k \nabla f(x^k)) \equiv F_{\gamma_k\nabla f}^h(x^k),
\ee
where $\gamma_k>0$ is the step size, and $h$ is a Legendre function. It was later showed by Beck and Teboulle \cite{BT-OR-letters-mirror-descent-2003} that the mirror descent algorithm \eqref{BGM-alg-explicit} can be interpreted as a projected gradient method with a Bregman distance, which can be described as follows:
\be\label{BGM-alg}
x^{k+1} := \argmin_{x\in \X} \ \langle \nabla f(x^k), x-x^k\rangle + \frac{1}{\gamma_k}D_h(x,x^k).
\ee
Here we discuss an important setting where $S$ is the probability simplex $S:=\{x\in\br^n\mid \sum_i x_i = 1, x\geq 0\}$. In this case, if one applies the projected gradient method using the Euclidean distance to solve \eqref{prob-uncon}, then each iteration requires a projection onto $S$. If one uses the Bregman distance generated by the Boltzmann-Shannon entropy, i.e., Example \ref{Example:Bregman-distance} (iii), then the solution of \eqref{BGM-alg} is given by a simple normalization
\[x^{k+1} := \frac{y^k}{\|y^k\|_1}, \quad \mbox{ with } y^k := x^k \circ e^{-\gamma_k\nabla f(x^k)}, \]
which is easier to compute than the projection onto $S$. Here for vectors $a$ and $b$, $a\circ b$ is the componentwise multipication, and $e^{a}$ is the componentwise exponential function.

\subsection{Bregman Proximal Gradient Method and Bregman Forward-Backward Splitting}
For composite convex minimization problem 
\be\label{prob-composite}
\min_{x\in S} \  f(x) + g(x) 
\ee
where $f$ is convex and smooth and $g$ is convex and possibly nonsmooth, Bregman proximal gradient method was studied in \cite{Teboulle-nolips,Lu-Bregman-PGM,Tseng-MPB-2010}, which iterates as: 
\begin{subequations}\label{BPGM-alg}
\begin{align}
x^{k+1}  & := \argmin_x \ \langle \nabla f(x^k), x-x^k\rangle + g(x) + \frac{1}{\gamma_k}D_h(x,x^k) \\
& \equiv \prox_{\gamma_kg}^h(\nabla h^*(\nabla h(x^k) - \gamma_k\nabla f(x^k))) \equiv J_{\gamma_k\partial g}^h \circ F_{\gamma_k\nabla f}^k (x^k),
\end{align}
\end{subequations}
where 
\be\label{def-Bregman-Prox-map}
\prox_g^h(z) := \argmin_x \ g(x) + D_h(x,z) \equiv J_{\partial g}^h(z),
\ee
is called the Bregman proximal map of $g$ with respect to $h$. An interesting question is how to accelerate the Bregman proximal gradient method using Nesterov's acceleration techniques \cite{Nesterov-1983,NesterovConvexBook2004,Nesterov-2005,Beck-Teboulle-2009}. When $f$ has an globally Lipschitz gradient, and $h$ is strongly convex, faster algorithms have been given by Auslender and Teboulle \cite{Auslender-Teboulle-2006}, and Tseng \cite{Tseng-MPB-2010}. However, when these assumptions are weakened, this problem has not been fully addressed in the literature. We refer to \cite{Hanzely-Richtarik-Xiao-2018,Gutman-Pena-Bregman-gradient} for some recent progresses on this topic and the recent paper by Teboulle \cite{Teboulle-survey-2018} for more detailed discussions on Bregman proximal gradient method. 

When it comes to the monotone inclusion problem \eqref{sum-2-operator} with $B$ being single-valued, the Bregman forward-backward splitting methods are studied in the literature, which iterates as 
\be\label{BFBS}
x^{k+1} := (\nabla h + \gamma_k A)^{-1}(\nabla h(x^k) - \gamma_k B(x^k)) \equiv J_{\gamma_k A}^h\circ F_{\gamma_k B}^h(x^k).
\ee
We refer to the recent paper by B\`ui and Combettes \cite{Bui-Combettes-Bregman-forward-backward} and references therein for more discussions on this method. 

\subsection{Bregman Proximal Point Method and Bregman Augmented Lagrangian method}
Another widely used Bregman algorithm is the Bregman proximal point method (PPM). The idea of Bregman PPM can be traced back to \cite{Erlander-1981,Eriksson-1985,Eggermont-1990}. The Bregman PPM in its current form was proposed by Censor and Zenios \cite{Censor-Zenios-1992} for convex minimization problem and by Eckstein \cite{Eckstein-93} for monotone inclusion problem. This method was further studied in \cite{Chen-Teboulle-93,Iusem-95,Ben-tal-Zibulevsky-1997,Teboulle-entropic-1992,Auslender-Teboulle-Ben-Tiba-MOR-1999,Iusem-Svaiter-Teboulle-1994,Iusem-Teboulle-1995,Bauschke-Borwein-Combettes-BregmanMonotoneOpt-2003}. For the maximal monotone inclusion problem $0\in T(x)$, the Bregman PPM iterates as
\be\label{Bregman-PPM-monotone-inclusion}
x^{k+1} := J_{\gamma_kT}^h(x^k) = (\nabla h + \gamma_k T)^{-1}\circ \nabla h(x^k).
\ee
For convex minimization problem \eqref{prob-uncon} with $f$ being nonsmooth, the Bregman PPM reduces to 
\be\label{Bregman-PPM-convex}
x^{k+1} := J_{\gamma_k\partial f}^h(x^k) = \argmin_x \ f(x) + \frac{1}{\gamma_k} D_h(x,x^k).
\ee
This algorithm generalizes the PPM in Euclidean space \cite{Martinet-70,Rockafellar-1976-ppa,Rockafellar1976,Iusem-ALM-survey-1999} to non-quadratic distance. Since it is usually difficult to solve the subproblem in \eqref{Bregman-PPM-monotone-inclusion} and \eqref{Bregman-PPM-convex} exactly, inexact Bregman PPM is studied in the literature \cite{Eckstein-iBPPA-1998,Solodov-Svaiter-2000,Yang-Toh-iBPPA-OT-2021}. Moreover, the exponential multiplier method proposed by Kort and Bertsekas \cite{Kort-Bertsekas-EMM-1972,Tseng-EMM} is known to be a special case of Bregman PPM with a specifically chosen $h$. We will discuss this method in more details in Section \ref{sec:ADEMM}. Furthermore, the nonlinear rescaling method developed by Polyak \cite{Polyak-1986-NR-Russian,Polyak-IBM-report-1989,Polyak-Griva-JOTA-2004,Griva-Polyak-MPA-2006} is also known to be equivalent to a Bregman PPM with suitably chosen distance and nonlinear penalty functions, as proved by Polyak and Teboulle \cite{Polyak-Teboulle-1997}. 

Similar to the connections between PPM and augmented Lagnrangain method (ALM) in the Euclidean case, there is a similar connection between Bregman PPM and Bregman ALM, as illustrated by Eckstein \cite{Eckstein-93}. Here we use the following convex optimization problem with linear equality constraint to illustrate the idea: 
\be\label{prob-convex-linear-equality}
\min_x  \ f(x), \ \st, \ Mx = b, x\in \X.
\ee
It can be shown that the Bregman PPM for solving the dual problem of \eqref{prob-convex-linear-equality} is equivalent to the following Bregman ALM for solving \eqref{prob-convex-linear-equality}:
\begin{subequations}\label{Bregman-ALM}
\begin{align}
x^{k+1} & \in \argmin_{x\in\X} \ f(x) + \frac{1}{\gamma_k}h^*(\nabla h(\lambda^k) + \gamma_k(Mx-b)) \\
\lambda^{k+1} & := \nabla h^*(\nabla h(\lambda^k)+\gamma_k(Mx^{k+1}-b)). 
\end{align}
\end{subequations}
One can see that unlike the classical ALM with a quadratic penalty function, the Bregman ALM \eqref{Bregman-ALM} adopts a non-quadratic penalty function $h^*$.
A classical reference on Bregman ALM is due to Bertsekas \cite{Bertsekas-book-96}, and a more recent survey is due to Iusem \cite{Iusem-ALM-survey-1999}. Some very recent works on this topic include \cite{Eckstein-practical-BALM,Yan-He-BALM}.

\subsection{Bregman ADMM}

A natural extension of the ADMM in Euclidean space is the Bregman ADMM. A Bregman ADMM was studied by Wang and Banerjee \cite{wang-bregman-admm}, in which the authors targeted solving \eqref{admm-sum-2} using the following Bregman ADMM: 
\begin{subequations}\label{wang-BADMM}
\begin{align}
u^{k+1} & := \argmin_y \ f(y) + \langle w^k, Mu+Nv^k-b\rangle + \gamma D_h(b-Mu, Nv^k)  \label{wang-BADMM-1}\\ 
v^{k+1} & := \argmin_z \ g(z) + \langle w^k, Mu^{k+1}+Nv-b\rangle + \gamma D_h(Nv,b-Mu^{k+1})  \label{wang-BADMM-2} \\
w^{k+1} & := w^k + \gamma(Mu^{k+1} + Nv^{k+1} - b). \label{wang-BADMM-3}
\end{align} 
\end{subequations}
Note that this algorithm requires $Nv^k$ and $b-Mu^{k+1}$ to lie in the domain of $h$, which is very difficult to guarantee for many useful Bregman distances such as the ones generated by the Boltzmann-Shannon entropy and the Burg's entropy. Thus the applicability of \eqref{wang-BADMM} is limited.

\section{Bregman Douglas-Rachford Splitting Method}\label{sec:BDRS}

In this section, we introduce our BDRS for solving \eqref{sum-2-operator}. A typical iteration of our BDRS for solving \eqref{sum-2-operator} is as follows: 
\begin{subequations}\label{BDRS-alg}
\begin{align}
x^k & := J_{\gamma_k B}^h (z^k) \\
y^k & := J_{\gamma_k A}^h\circ \nabla h^*(2\nabla h(x^k)-\nabla h(z^k)) \\
z^{k+1} & := \nabla h^*(\nabla h(z^k) - \nabla h(x^k) + \nabla h(y^k)),
\end{align}
\end{subequations} 
where $\gamma_k>0$ is a parameter. 
We now provide some explainations to our BDRS \eqref{BDRS-alg}. Note that the Bregman resolvent operators are always taken to the points in the primal space. For points in the mirror space, we always use $\nabla h^*$ to transform them back to the primal space. Note that \eqref{BDRS-alg} can be written equivalently as 
\be\label{BDRS-alg-compact}
\nabla h(z^{k+1}) = \nabla h(z^k) + \nabla h\circ J_{\gamma_k A}^h\circ \nabla h^*(2\nabla h\circ J_{\gamma_k B}^h(z^k)-\nabla h(z^k)) - \nabla h\circ J_{\gamma_k B}^h(z^k).
\ee
Using the Bregman reflection operator \eqref{def-reflection-op} and the Bregman Mann's operator \eqref{def-Mann}, the BDRS \eqref{BDRS-alg-compact} can be  written more compactly as
\be\label{BDRS-alg-Mann}
z^{k+1} := M_{\half}^h(R_{\gamma_k A}^hR_{\gamma_k B}^h) (z^k).
\ee
We notice that when $h(\cdot) = \half\|\cdot\|_2^2$, all the three forms of BDRS, i.e., \eqref{BDRS-alg}, \eqref{BDRS-alg-compact}, and \eqref{BDRS-alg-Mann} reduce exactly to their Euclidean counterparts.

\subsection{Application to Convex Minimization: Bregman ADMM}

First, we note that if one applies the Bregman ALM to solve the convex optimization problem with linear equality constraints \eqref{admm-sum-2}, then it should iterate as follows:
\begin{subequations}\label{Bregman-ALM-2-block}
\begin{align}
(u^k,v^k) & := \argmin_{u,v} \ f(u) + g(v)+ \frac{1}{\gamma_k}h^*(\nabla h(w^k) + \gamma_k(Mu+Nv-b)) \\
w^{k+1} & := \nabla h^*(\nabla h(w^k)+\gamma_k(Mu^k+Nv^k-b)). 
\end{align}
\end{subequations}
As a result, it is easy to see that the Bregman ADMM for solving \eqref{admm-sum-2} is given by the following updates: 
\begin{subequations}\label{Bregman-ADMM}
\begin{align}
u^k & := \argmin_{u} \ f(u) + \frac{1}{\gamma_k}h^*(\nabla h(w^k) + \gamma_k(Mu+Nv^{k-1}-b)) \\
v^k & := \argmin_{v} \ g(v) + \frac{1}{\gamma_k}h^*(\nabla h(w^k) + \gamma_k(Mu^k+Nv-b)) \\
w^{k+1} & := \nabla h^*(\nabla h(w^k)+\gamma_k(Mu^k+Nv^k-b)),
\end{align}
\end{subequations}
where we alternatingly update $u^k$ and $v^k$ in Bregman ALM with the other variable being fixed. 
To the best of our knowledge, both BDRS \eqref{BDRS-alg} and Bregman ADMM \eqref{Bregman-ADMM} are new to the literature. In the following, we show that the Bregman ADMM \eqref{Bregman-ADMM} is actually a direct application of BDRS \eqref{BDRS-alg} to the dual of \eqref{admm-sum-2}, which is given by 
\begin{subequations}\label{admm-sum-2-dual}
\begin{align}
\max_w\min_{u,v} f(u) + g(v) + \langle w, Mu+Nv-b\rangle & \equiv \max_w -f^*(-M^\top w) - g^*(-N^\top w) - \langle b,w\rangle \\
& \equiv \min_w f^*(-M^\top w) + g^*(-N^\top w) + \langle b,w\rangle,
\end{align}
\end{subequations}
where $w$ is the dual variable (Lagrange multiplier), and $f^*$ and $g^*$ are the conjugate functions of $f$ and $g$, respectively. The optimality condition of \eqref{admm-sum-2-dual} is given by 
\be\label{admm-sum-2-dual-OptCond}
0 \in -M \partial f^*(-M^\top w) - N\partial g^*(-N^\top w) + b,
\ee
which is in the from of \eqref{sum-2-operator} with $A(w) = -M \partial f^*(-M^\top w)$ and $B(w) = -N\partial g^*(-N^\top w) + b$. 

\begin{theorem}\label{thm:BDRS-equiv-BADMM-linear-equality}
The BDRS \eqref{BDRS-alg} for solving the dual problem \eqref{admm-sum-2-dual-OptCond} is equivalent to the Bregman ADMM \eqref{Bregman-ADMM} for solving the primal problem \eqref{admm-sum-2}.
\end{theorem}
\begin{proof}
With $A(w) = - M \partial f^*(-M^\top w)$ and $B(w) = - N\partial g^*(-N^\top w) +  b$, the BDRS for solving \eqref{admm-sum-2-dual-OptCond} can be written as
\begin{subequations}\label{BDRS-admm-sum-2-dual-OptCond}
\begin{align}
x^k & := \argmin_x \ g^*(-N^\top x) + b^\top x + \frac{1}{\gamma_k} D_h(x,z^k) \label{BDRS-admm-sum-2-dual-OptCond-1}\\
y^k & := \argmin_y \ f^*(-M^\top y) +\frac{1}{\gamma_k} D_h(y, \nabla h^*(2\nabla h(x^k)-\nabla h(z^k))) \label{BDRS-admm-sum-2-dual-OptCond-2}\\
z^{k+1} & := \nabla h^*(\nabla h(z^k) - \nabla h(x^k) + \nabla h(y^k)). \label{BDRS-admm-sum-2-dual-OptCond-3}
\end{align}
\end{subequations}
By using $v$ to denote the dual variable for \eqref{BDRS-admm-sum-2-dual-OptCond-1}, we obtain that \eqref{BDRS-admm-sum-2-dual-OptCond-1} is equivalent to 
\begin{subequations}\label{BDRS-admm-sum-2-dual-OptCond-1-equiv}
\begin{align}
    & \min_x \max_v \ \langle -N^\top x, v\rangle - g(v) + b^\top x + \frac{1}{\gamma_k}D_h(x,z^k) \\
= & \max_v \min_x \ \langle -N^\top x, v\rangle - g(v) + b^\top x + \frac{1}{\gamma_k}D_h(x,z^k) \\
\equiv & \max_v \min_x \ h(x) - \langle \nabla h(z^k)+\gamma_k(Nv-b), x\rangle -\gamma_k g(v) - h(z^k) + \langle \nabla h(z^k), z^k \rangle \\
= & \max_v -\{\max_x -h(x) + \langle \nabla h(z^k)+\gamma_k(Nv-b), x\rangle +\gamma_k g(v) + h(z^k) - \langle \nabla h(z^k), z^k \rangle\} \label{BDRS-admm-sum-2-dual-OptCond-1-equiv-4}\\
= & \max_v -\{h^*(\nabla h(z^k)+\gamma_k(Nv-b)) +\gamma_k g(v) + h(z^k) - \langle \nabla h(z^k), z^k \rangle\} \\
\equiv & \min_v h^*(\nabla h(z^k)+\gamma_k(Nv-b)) +\gamma_k g(v).
\end{align}
\end{subequations}
Moreover, note that the optimal $x$ in \eqref{BDRS-admm-sum-2-dual-OptCond-1-equiv-4} is given by:
\be\label{BDRS-admm-sum-2-dual-OptCond-1-equiv-4-OptCond}
x := \nabla h^*(\gamma_k(Nv-b)+\nabla h(z^k)).
\ee
Similarly, by using $u$ to denote the dual variable for \eqref{BDRS-admm-sum-2-dual-OptCond-2}, we obtain that \eqref{BDRS-admm-sum-2-dual-OptCond-2} is equivalent to (for ease of presentation, denote $\tilde{y}^k:=\nabla h^*(2\nabla h(x^k)-\nabla h(z^k))$): 
\begin{subequations}\label{BDRS-admm-sum-2-dual-OptCond-2-equiv}
\begin{align}
    & \min_y \max_u \ \langle -M^\top y, u\rangle - f(u) + \frac{1}{\gamma_k}D_h(y,\tilde{y}^k) \\
= & \max_u \min_y \ \langle -M^\top y, u\rangle - f(u) + \frac{1}{\gamma_k} D_h(y,\tilde{y}^k) \\
\equiv & \max_u \min_y \ h(y) - \langle \nabla h(\tilde{y}^k)+\gamma_k Mu, y\rangle -\gamma_k f(u) - h(\tilde{y}^k) + \langle \nabla h(\tilde{y}^k), \tilde{y}^k \rangle \\
= & \max_u -\{\max_y -h(y) + \langle \nabla h(\tilde{y}^k)+\gamma_k Mu, y\rangle +\gamma_k f(u) + h(\tilde{y}^k) - \langle \nabla h(\tilde{y}^k), \tilde{y}^k \rangle\} \label{BDRS-admm-sum-2-dual-OptCond-2-equiv-4}\\
= & \max_u -\{h^*(\nabla h(\tilde{y}^k)+\gamma_k Mu) +\gamma_k f(u) + h(\tilde{y}^k) - \langle \nabla h(\tilde{y}^k), \tilde{y}^k \rangle\} \\
\equiv & \min_u h^*(2\nabla h(x^k)-\nabla h(z^k)+\gamma_k Mu) +\gamma_k f(u).
\end{align}
\end{subequations}
Moreover, the optimal $y$ in \eqref{BDRS-admm-sum-2-dual-OptCond-2-equiv-4} is given by
\be\label{BDRS-admm-sum-2-dual-OptCond-2-equiv-4-OptCond}
y := \nabla h^*(2\nabla h(x^k)-\nabla h(z^k)+\gamma_k Mu).
\ee
By combining \eqref{BDRS-admm-sum-2-dual-OptCond-1-equiv}, \eqref{BDRS-admm-sum-2-dual-OptCond-1-equiv-4-OptCond}, \eqref{BDRS-admm-sum-2-dual-OptCond-2-equiv} and \eqref{BDRS-admm-sum-2-dual-OptCond-2-equiv-4-OptCond}, we have that the BDRS \eqref{BDRS-admm-sum-2-dual-OptCond} is equivalent to
\begin{subequations}\label{BDRS-admm-sum-2-dual-OptCond-rewrite-uv}
\begin{align}
v^k & := \argmin_v h^*(\nabla h(z^k)+\gamma_k(Nv-b)) +\gamma_k g(v) \\
x^k & := \nabla h^*(\gamma_k(Nv^k-b)+\nabla h(z^k)) \label{BDRS-admm-sum-2-dual-OptCond-rewrite-uv-x}\\
u^k & := \argmin_u h^*(2\nabla h(x^k)-\nabla h(z^k)+\gamma_k Mu) +\gamma_k f(u) \\
y^k & := \nabla h^*(2\nabla h(x^k)-\nabla h(z^k)+\gamma_k Mu^k) \label{BDRS-admm-sum-2-dual-OptCond-rewrite-uv-y}\\
z^{k+1} & := \nabla h^*(\nabla h(z^k) - \nabla h(x^k) + \nabla h(y^k)).\label{BDRS-admm-sum-2-dual-OptCond-rewrite-uv-z}
\end{align}
\end{subequations}
Reordering the updates in \eqref{BDRS-admm-sum-2-dual-OptCond-rewrite-uv}, we know that \eqref{BDRS-admm-sum-2-dual-OptCond-rewrite-uv} is equivalent to 
\begin{subequations}\label{BDRS-admm-sum-2-dual-OptCond-rewrite-uv-reorder}
\begin{align}
u^k & := \argmin_u h^*(2\nabla h(x^k)-\nabla h(z^k)+\gamma_k Mu) +\gamma_k f(u) \label{BDRS-admm-sum-2-dual-OptCond-rewrite-uv-reorder-1}\\
y^k & := \nabla h^*(2\nabla h(x^k)-\nabla h(z^k)+\gamma_k Mu^k) \label{BDRS-admm-sum-2-dual-OptCond-rewrite-uv-reorder-2}\\
z^{k+1} & := \nabla h^*(\nabla h(z^k) - \nabla h(x^k) + \nabla h(y^k)) \label{BDRS-admm-sum-2-dual-OptCond-rewrite-uv-reorder-3}\\
v^{k+1} & := \argmin_v h^*(\nabla h(z^{k+1})+\gamma_k(Nv-b)) +\gamma_k g(v) \label{BDRS-admm-sum-2-dual-OptCond-rewrite-uv-reorder-4}\\
x^{k+1} & := \nabla h^*(\gamma_k(Nv^{k+1}-b)+\nabla h(z^{k+1})). \label{BDRS-admm-sum-2-dual-OptCond-rewrite-uv-reorder-5}
\end{align}
\end{subequations}
Combining \eqref{BDRS-admm-sum-2-dual-OptCond-rewrite-uv-reorder-2} and \eqref{BDRS-admm-sum-2-dual-OptCond-rewrite-uv-reorder-3} yields 
\be\label{reorder-relation-1}
\nabla h(z^{k+1}) = \nabla h(z^k) - \nabla (x^k) + 2\nabla h(x^k) - \nabla h(z^k) + \gamma_k Mu^k = \nabla h(x^k) +\gamma_k Mu^k.
\ee
From \eqref{BDRS-admm-sum-2-dual-OptCond-rewrite-uv-reorder-5} we have
\be\label{reorder-relation-2}
\nabla h(z^{k+1}) = \nabla h(x^{k+1}) - \gamma_k (Nv^{k+1}-b),
\ee
which implies 
\be\label{reorder-relation-3}
2\nabla h(x^{k+1}) - \nabla h(z^{k+1}) = \nabla h(x^{k+1}) +\gamma_k (Nv^{k+1}-b),
\ee
and 
\be\label{reorder-relation-4}
\nabla h(x^{k+1}) = \nabla h(x^k) +\gamma_k (Mu^k+Nv^{k+1}-b).
\ee
By substituting \eqref{reorder-relation-3} to \eqref{BDRS-admm-sum-2-dual-OptCond-rewrite-uv-reorder-1}, \eqref{reorder-relation-1} to \eqref{BDRS-admm-sum-2-dual-OptCond-rewrite-uv-reorder-4}, and combining with \eqref{reorder-relation-4}, we know that \eqref{BDRS-admm-sum-2-dual-OptCond-rewrite-uv-reorder} is equivalent to
\begin{subequations}\label{BDRS-admm-sum-2-dual-OptCond-rewrite-uv-reorder-simplified}
\begin{align}
u^k & := \argmin_u h^*(\nabla h(x^k) + \gamma_k(Nv^k-b+Mu)) +\gamma_k f(u)  \\ 
v^{k+1} & := \argmin_v h^*(\nabla h(x^k) + \gamma_k(Mu^k+Nv-b)) +\gamma_k g(v) \\
x^{k+1} & := \nabla h^*(\nabla h(x^k) + \gamma_k(Mu^k+Nv^{k+1}-b)).
\end{align}
\end{subequations}
This is exactly the same as the Bregman ADMM \eqref{Bregman-ADMM}.
\end{proof}

\subsection{Connection to Variable Metric ADMM}

A variable metric ADMM has been studied in \cite{Giselsson-variable-metric-DRS-ADMM,Bot-variable-metric-ADMM}, which solves the convex minimization problem \eqref{admm-sum-2} using the following updates:
\begin{subequations}\label{BDRS-alg-L-equiv-dual-switch-uv}
\begin{align}
u^k & := \argmin_u \ f(u) + \frac{1}{2\gamma_k}\|\gamma_k(Mu + Nv^{k-1}-b) + Lw^k\|_{L^{-1}}^2 \\
v^k & := \argmin_v \ g(v) + \frac{1}{2\gamma_k}\|\gamma_k(Mu^k+Nv-b) +Lw^k\|_{L^{-1}}^2 \\
w^{k+1} & := w^k + \gamma_kL^{-1}(Mu^k + Nv^k-b),
\end{align}
\end{subequations}
where $L \succ 0$ is a positive definite matrix, and $\|x\|_L^2 := x^\top Lx$. This algorithm is also known as applying ADMM to \eqref{admm-sum-2} with preconditioned constraints \cite{Ghadimi-variable-metric-admm}. 
We now show that this variable metric ADMM is a special case of our Bregman ADMM \eqref{Bregman-ADMM} with a special choice of $h(x) = \half\|x\|_L^2 := \half x^\top L x$, i.e., the quadratic function in Example \ref{Example:Bregman-distance} (ii). In this case, we have
\be\label{VM-ADMM-h}
\nabla h(x) = Lx, \quad h^*(x) = \half\|x\|_{L^{-1}}^2 = \half x^\top L^{-1}x, \quad \nabla h^*(x) = L^{-1}x. 
\ee

\begin{theorem}\label{thm:VM-ADMM}
The variable metric ADMM \eqref{BDRS-alg-L-equiv-dual-switch-uv} is a special case of the Bregman ADMM \eqref{Bregman-ADMM} with $h$ given in Example \ref{Example:Bregman-distance} (ii).
\end{theorem}

\begin{proof}
From \eqref{VM-ADMM-h}, we know that the Bregman ADMM \eqref{Bregman-ADMM} can be written as 
\begin{subequations}\label{Bregman-ADMM-VM}
\begin{align}
u^k & := \argmin_{u} \ f(u) + \frac{1}{2\gamma_k}\|Lw^k + \gamma_k(Mu+Nv^{k-1}-b)\|_{L^{-1}}^2 \\
v^k & := \argmin_{v} \ g(v) + \frac{1}{2\gamma_k}\|Lw^k + \gamma_k(Mu^k+Nv-b)\|_{L^{-1}}^2 \\
w^{k+1} & := L^{-1}(Lw^k+\gamma_k(Mu^k+Nv^k-b)),
\end{align}
\end{subequations}
which is the same as the variable metric ADMM \eqref{BDRS-alg-L-equiv-dual-switch-uv}.
\end{proof}

\subsection{Alternating Direction Exponential Multiplier Method}\label{sec:ADEMM}

In this section, we propose a new algorithm, which is a special case of BDRS, for solving the following linear inequality constrained convex minimization problem: 
\be\label{admm-sum-2-inequality-constraint}
\min_{u,v} \ f(u) + g(v), \ \st, \ M u+ N v-b \leq 0, \ u\in\br^p, v\in\br^q,
\ee
where $f$ and $g$ are both proper, closed and convex functions, and $M\in\br^{m\times p}$, $N\in\br^{m\times q}$, $b\in \br^m$.

One important approach for solving \eqref{admm-sum-2-inequality-constraint} is the exponential multiplier method (EMM) that was proposed and studied by Bertsekas, Kort and Tseng \cite{Kort-Bertsekas-EMM-1972,Bertsekas-book-96,Tseng-EMM}. Unlike the usual augmented Lagrangian method that uses a quadratic penalty term, the EMM proposes to use a non-quadratic penalty term.
More specifically, the EMM uses the exponential penalty function given by:
\be\label{def-psi}
\psi(t) = e^t - 1.
\ee
By associating a Lagrange multiplier $w_j$ to the $j$-th constraint in \eqref{admm-sum-2-inequality-constraint}, the EMM for solving \eqref{admm-sum-2-inequality-constraint} iterates as follows \footnote{Note that the EMM proposed in \cite{Kort-Bertsekas-EMM-1972,Bertsekas-book-96,Tseng-EMM} allows the scalar $\gamma_k$ to be replaced by a vector whose $j$-th entry is $[\gamma_k]_j$. In this case, the EMM becomes \begin{align*}
(u^k,v^k) & := \argmin_{u,v} f(u) + g(v) + \sum_{j=1}^m \frac{w_j^k}{\gamma_j^k}\psi([\gamma_k]_j(M_j^\top u+N_j^\top v-b_j)) \\ 
w_j^{k+1} & := w_j^k\nabla\psi([\gamma_k]_j(M_j^\top u^k+N_j^\top v^k-b_j)) = w_j^k e^{[\gamma_k]_j(M_j^\top u^k+N_j^\top v^k-b_j)}, \ j=1,\ldots,m.
\end{align*} Here we use a scalar $\gamma_k$ for simplicity.}:
 
\begin{subequations}\label{EMM-alg}
\begin{align}
(u^k,v^k) & := \argmin_{u,v} f(u) + g(v) + \frac{1}{\gamma_k}\sum_{j=1}^m w_j^k\psi(\gamma_k(M_j^\top u+N_j^\top v-b_j)) \label{EMM-alg-1} \\ 
w_j^{k+1} & := w_j^k\nabla\psi(\gamma_k(M_j^\top u^k+N_j^\top v^k-b_j)) = w_j^k e^{\gamma_k(M_j^\top u^k+N_j^\top v^k-b_j)}, \ j=1,\ldots,m,\label{EMM-alg-2}
\end{align}
\end{subequations}
where $\gamma_k > 0$ is a parameter, and $M_j^\top$ denotes the $j$-th row of $M$ and $N_j^\top$ denotes the $j$-th row of $N$. It can be shown that the EMM \eqref{EMM-alg} is equivalent to a Bregman proximal point algorithm for solving the dual of \eqref{admm-sum-2-inequality-constraint}. Note that the dual of \eqref{admm-sum-2-inequality-constraint} is given by
\be\label{admm-sum-2-inequality-constraint-dual}
\max \ d(w), \ \st, \ w\geq 0,
\ee
where the dual function $d(w)$ is defined as
\[
d(w) := \min_{u,v} \left\{f(u)+g(v) + \sum_{j=1}^m w_j (M_j^\top u+N_j^\top v-b_j)\right\}.
\]
It can be shown that the EMM \eqref{EMM-alg} for solving the primal problem \eqref{admm-sum-2-inequality-constraint} is equivalent to the following algorithm for solving the dual problem \eqref{admm-sum-2-inequality-constraint-dual}:
\be\label{EMM-alg-dual-BPPA}
w^{k+1} := \argmax_{w\geq 0} \ \left\{d(w) - \frac{1}{\gamma_k}\sum_{j=1}^m w_j^k\psi^*\left(\frac{w_j}{w_j^k}\right)\right\},
\ee
where $\psi^*$ is the conjugate function of $\psi$, which is the entropy function 
\be\label{def-psi-conj}
\psi^*(s) = s\log s - s + 1. 
\ee
Note that \eqref{EMM-alg-dual-BPPA} is indeed a Bregman proximal point algorithm, because $\sum_{j=1}^m{w_j^k}\psi^*\left(\frac{w_j}{w_j^k}\right) = D_h(w,w^k)$ with $h$ being the Boltzmann-Shannon entropy defined in Example \eqref{Example:Bregman-distance} (iii).

Note that for some applications, the subproblem \eqref{EMM-alg-1} in EMM can be difficult to solve. To overcome this difficulty, we propose an alternating direction exponential multiplier method (ADEMM) for solving \eqref{admm-sum-2-inequality-constraint}. The ADEMM iterates as follows: 
\begin{subequations}\label{ADEMM-alg}
\begin{align}
u^k & := \argmin_{u} f(u) + \frac{1}{\gamma_k}\sum_{j=1}^m w_j^k\psi(\gamma_k(M_j^\top u+N_j^\top v^{k-1}-b_j)) \label{ADEMM-alg-1} \\ 
v^k & := \argmin_{v} g(v) + \frac{1}{\gamma_k}\sum_{j=1}^m w_j^k\psi(\gamma_k(M_j^\top u^k+N_j^\top v-b_j)) \label{ADEMM-alg-2} \\ 
w_j^{k+1} & := w_j^k\nabla\psi(\gamma_k(M_j^\top u^k+N_j^\top v^k-b_j)) = w_j^k e^{\gamma_k(M_j^\top u^k+N_j^\top v^k-b_j)}, \ j=1,\ldots,m.\label{ADEMM-alg-3}
\end{align}
\end{subequations}
That is, the ADEMM alternatingly minimizes the Lagrangian function with respect to $u$ and $v$ in each iteration. This is exactly in the same spirit of the usual ADMM in Euclidean space with a quadratic penalty. To the best of our knowledge, the ADEMM \eqref{ADEMM-alg} is new to the literature: it is the first alternating direction version of EMM. For some applications, both subproblems \eqref{ADEMM-alg-1} and \eqref{ADEMM-alg-2} are easier to solve than  \eqref{EMM-alg-1}. This is indeed the case as we will see later in the discrete optimal transport problem in Section \ref{sec:OT}. 

Very interestingly, the ADEMM \eqref{ADEMM-alg} for solving the primal problem  \eqref{admm-sum-2-inequality-constraint} is equivalent to BDRS \eqref{BDRS-alg} for solving the dual problem \eqref{admm-sum-2-inequality-constraint-dual}, with $h$ being the Boltzmann-Shannon entropy.

\begin{theorem}
The ADEMM \eqref{ADEMM-alg} for solving the primal problem \eqref{admm-sum-2-inequality-constraint} is equivalent to BDRS \eqref{BDRS-alg} for solving the dual problem \eqref{admm-sum-2-inequality-constraint-dual}, with $h$ being the Boltzmann-Shannon entropy.
\end{theorem}

\begin{proof}
We first note that the dual problem \eqref{admm-sum-2-inequality-constraint-dual} is equivalent to 
\be\label{admm-sum-2-inequality-constraint-dual-equivalent}
\min_{w\geq 0} \ f^*(-M^\top w) + g^*(-N^\top w) + b^\top w,
\ee
whose optimality condition is given by:
\be\label{admm-sum-2-inequality-constraint-dual-equivalent-OptCond}
0 \in -M \partial f^*(M^\top w) - N \partial g^*(N^\top w) + b + \partial \bone(w\geq 0),
\ee
where $\bone(C)$ denotes the indicator function of set $C$. By letting 
$A(w) = - M \partial f^*(-M^\top w) + \partial \bone(w\geq 0)$ and $B(w) = - N\partial g^*(-N^\top w) + b + \partial \bone(w\geq 0)$, we note that \eqref{admm-sum-2-inequality-constraint-dual-equivalent-OptCond} is in the same form as \eqref{sum-2-operator} and thus can be solved by BDRS \eqref{BDRS-alg}.
The BDRS \eqref{BDRS-alg} with $h$ being the Boltzmann-Shannon entropy can be written as
\begin{subequations}\label{BDRS-alg-ADEMM}
\begin{align}
x^k & := \argmin_{x\geq 0} \ g^*(-N^\top x) + b^\top x + \frac{1}{\gamma_k} D_h(x,z^k) \label{BDRS-alg-ADEMM-1} \\
y^k & := \argmin_{y\geq 0} \ f^*(-M^\top y) +\frac{1}{\gamma_k} D_h(y, (x^k\circ x^k)/z^k) \label{BDRS-alg-ADEMM-2}\\
z^{k+1} & := \nabla h^*(\nabla h(z^k) - \nabla h(x^k) + \nabla h(y^k)),\label{BDRS-alg-ADEMM-3}
\end{align}
\end{subequations}
where $a\circ b$ denotes the elementwise multiplication of vectors $a$ and $b$, and $a/b$ denotes the elementwise division of vectors $a$ and $b$. 

By introducing $v$ as the dual variable of \eqref{BDRS-alg-ADEMM-1}, we obtain that \eqref{BDRS-alg-ADEMM-1} is equivalent to:
\begin{align*}
     & \min_x \max_v \langle -N^\top x, v\rangle - g(v) + b^\top x +\frac{1}{\gamma_k} D_h(x,z^k) \\
=  & \max_v \min_x \langle -N^\top x, v\rangle - g(v) + b^\top x + \frac{1}{\gamma_k} D_h(x,z^k) \\
\equiv & \max_v -\gamma_k g(v) - h^*(\nabla h(z^k) + \gamma_k (N_j^\top v-b_j)) \\
\equiv & \min_v g(v) + \frac{1}{\gamma_k}\sum_j z^k_j e^{\gamma_k (N_j^\top v-b_j)},
\end{align*}
with $\gamma_k(-Nv + b) + \nabla h(x) - \nabla h(z^k)=0$. Similarly, by introducing $u$ as the dual variable of \eqref{BDRS-alg-ADEMM-2}, we obtain that \eqref{BDRS-alg-ADEMM-2} is equivalent to:
\begin{align*}
     & \min_y \max_u \langle -M^\top y, u\rangle - f(u) + \frac{1}{\gamma_k}D_h(y, (x^k\circ x^k)/z^k) \\
=  & \max_u \min_y \langle -M^\top y, u\rangle - f(u) + \frac{1}{\gamma_k} D_h(y, (x^k\circ x^k)/z^k) \\
=  & \max_u -\gamma f(u) - \sum_j  (x^k_j \cdot x^k_j / z^k_j) e^{\gamma_k M_j^\top u} \\
\equiv & \min_u f(u) + \frac{1}{\gamma_k} \sum_j  (x^k_j \cdot x^k_j / z^k_j) e^{\gamma_k M_j^\top u},
\end{align*}
with $-\gamma_k Mu+\nabla h(y) - \nabla h((x^k\circ x^k)./z^k)=0$. Therefore, \eqref{BDRS-alg-ADEMM} can be equivalently rewritten as
\begin{subequations}\label{BDRS-alg-ADEMM-equivalent}
\begin{align}
v^k & := \argmin_v g(v) + \frac{1}{\gamma_k}\sum_j z^k_j e^{\gamma_k(N_j^\top v-b_j)} \label{BDRS-alg-ADEMM-equivalent-1} \\
x^k_j & := z^k_j e^{\gamma_k(N_j^\top v^k-b_j)} \\ 
u^k & := \argmin_u f(u) + \frac{1}{\gamma_k}\sum_j  (x^k_j \cdot x^k_j / z^k_j) e^{\gamma_k M_j^\top u} \label{BDRS-alg-ADEMM-equivalent-2}\\
y^k_j & := (x^k_j \cdot x^k_j / z^k_j) e^{\gamma_k M_j^\top u^k}\\
z^{k+1} & := \nabla h^*(\nabla h(z^k) - \nabla h(x^k) + \nabla h(y^k)).\label{BDRS-alg-ADEMM-equivalent-3}
\end{align}
\end{subequations}
We have \eqref{BDRS-alg-ADEMM-equivalent} can be equivalently rewritten as 
\begin{align*}
v^k & := \argmin_v g(v) + \frac{1}{\gamma_k}\sum_j z_j^k e^{\gamma_k(N_j^\top v-b_j)}\\
x^k_j & := z^k_j e^{\gamma_k(N_j^\top v^k-b_j)} \\ 
u^k & := \argmin_u f(u) + \frac{1}{\gamma_k}\sum_j  x^k_j e^{\gamma_k(M_j^\top u+N_j^\top v^k-b_j)}\\ 
z^{k+1}_j & := \nabla h^*(\nabla h(z^k_j) + \gamma_k(M_j^\top u^k+N_j^\top v^k-b_j)).
\end{align*}
The last equation implies that $\nabla h(z_j^{k+1}) =\nabla h(x_j^k) + \gamma_k M_j^\top u^k$, and therefore, $z_j^{k+1} = x_j^k e^{\gamma_k M_j^\top u^k}$. Thus, the above is equivalent to
\begin{align*}
v^k & := \argmin_v g(v) + \frac{1}{\gamma_k} \sum_j x_j^{k-1} e^{\gamma_k(M_j^\top u^{k-1}+N_j^\top v-b_j)}\\
x^k_j & := x^{k-1}_j e^{\gamma_k(M_j^\top u^{k-1}+N_j^\top v^k-b_j)} \\ 
u^k & := \argmin_u f(u) +  \frac{1}{\gamma_k} \sum_j  x^k_j e^{\gamma_k(M_j^\top u+N_j^\top v^k-b_j)} .
\end{align*}
This is exactly the ADEMM \eqref{ADEMM-alg}. 
\end{proof}

We end this section by remarking that the EMM \eqref{EMM-alg}, being equivalent to a Bregman PPM, is also equivalent to the nonlinear rescaling method developed by Polyak \cite{Polyak-1986-NR-Russian,Polyak-IBM-report-1989,Polyak-Griva-JOTA-2004,Griva-Polyak-MPA-2006} with suitably chosen nonlinear rescaling function, as proved by Polyak and Teboulle \cite{Polyak-Teboulle-1997}. More specifically, the linearly inequality constrained problem \eqref{admm-sum-2-inequality-constraint} is equivalent to:
\be\label{admm-sum-2-inequality-constraint-NR}
\min_{u,v} \ f(u) + g(v), \ \st, \ \frac{1}{\gamma_k}\psi(\gamma_k(M_j^\top u+N_j^\top v-b_j)) \leq 0, j=1,\ldots,m,
\ee
where $\psi$ is defined in \eqref{def-psi}. By associating a Lagrange multiplier $w_j$ to the $j$-th constraint in \eqref{admm-sum-2-inequality-constraint-NR}, the Lagrangian function of \eqref{admm-sum-2-inequality-constraint-NR} is given by:
\be\label{admm-sum-2-inequality-constraint-NR-Lag-Fn}
\LCal_{NR}(u,v;w) := f(u) + g(v) + \frac{1}{\gamma_k}\sum_j w_j\psi(\gamma_k(M_j^\top u+N_j^\top v-b_j)).
\ee
The nonlinear rescaling method is essentially a Lagrangian multiplier method for solving \eqref{admm-sum-2-inequality-constraint-NR-Lag-Fn} and is given by:
\begin{subequations}\label{admm-sum-2-inequality-constraint-NR-ALM}
\begin{align}
(u^k,v^k) & := \argmin_{u,v} \ \LCal_{NR}(u,v;w^k) \label{admm-sum-2-inequality-constraint-NR-ALM-1}\\
w_j^{k+1} & := w_j^k \nabla \psi(\gamma_k(M_j^\top u^k+N_j^\top v^k-b_j)), j=1,\ldots,m.\label{admm-sum-2-inequality-constraint-NR-ALM-2}
\end{align}
\end{subequations}
Polyak and Griva \cite{Polyak-Griva-JOTA-2004,Griva-Polyak-MPA-2006} proposed to use Newton's method to solve a primal-dual system that consists \eqref{admm-sum-2-inequality-constraint-NR-ALM-2} and the KKT system of \eqref{admm-sum-2-inequality-constraint-NR-ALM-1}. Note that Newton's method can be employed because function $\psi$ and its derivative $\psi'$ are both twice continuously differentiable. This is one of the main motivations of designing the nonlinear rescaling method. Our ADEMM \eqref{ADEMM-alg} leads to the following alternating direction version of the nonlinear rescaling method:
\begin{align*}
u^k & := \argmin_{u,v} \ \LCal_{NR}(u,v^{k-1};w^k) \\
v^k & := \argmin_{u,v} \ \LCal_{NR}(u^k,v;w^k) \\
w_j^{k+1} & := w_j^k \nabla \psi(\gamma_k(M_j^\top u^k+N_j^\top v^k-b_j)), j=1,\ldots,m. 
\end{align*}

\section{Bregman Peaceman-Rachford Splitting and Bregman Double-Backward Method}\label{sec:BPRS}

In this section, we discuss two algorithms that are related to BDRS, namely Bregman Peaceman-Rachford spliting method (BPRS) and Bregman double-backward method (BDBM). The Peaceman-Rachford splitting (PRS) method is another well-studied operator splitting method that was also proposed to solve variational problems arising from numerical PDEs \cite{Peaceman-Rachford-55}. The PRS can also be applied to solving the monotone inclusion problem \eqref{sum-2-operator}. Using our notions defined in Definition \ref{def-Bregman-operators}, the BPRS for solving \eqref{sum-2-operator} can be written as 
\be\label{BPRS-alg}
z^{k+1} := R_{\gamma_k A}^hR_{\gamma_k B}^h(z^k).
\ee
When $h(x)=\half\|x\|_2^2$, \eqref{BPRS-alg} reduces to the original PRS in the Euclidean space. For the convex minimization problem \eqref{admm-sum-2}, it is known that a symmetric ADMM is equivalent to the PRS applied to solving the dual of \eqref{admm-sum-2} \cite{Glowinski-LeTallec-89}. The symmetric ADMM for solving \eqref{admm-sum-2} iterates as follows: 
\begin{subequations}\label{sym-admm-sum-2-alg}
\begin{align}
u^{k+1} & := \argmin_u \ \LCal_\beta(u,v^k;w^k) \\
w^{k+\half} & := w^k + \beta(Mu^{k+1} + Nv^{k} - b)\\
v^{k+1} & := \argmin_v \ \LCal_\beta(u^{k+1},v;w^{k+\half})\\
w^{k+1} & := w^{k+\half} + \beta(Mu^{k+1} + Nv^{k+1} - b),
\end{align}
\end{subequations}
where the augmented Lagrangian function $\LCal_\beta$ is defined in \eqref{admm-sum-2-AL-Fn}. This symmetric ADMM is equivalent to an alternating linearization method when $f$ and $g$ are both smooth, as studied by Goldfarb, Ma and Scheinberg \cite{Goldfarb-Ma-Scheinberg-2010}. For BPRS \eqref{BPRS-alg}, we can prove a similar result. 

\begin{theorem}\label{thm:BPRS-equiv-symBADMM}
For $\gamma_k>0$, the BPRS \eqref{BPRS-alg} for solving the dual of \eqref{admm-sum-2} with $A(x) = -M\partial f^*(-M^\top x)$ and $B(x) = -N \partial g^*(-N^\top x) + b$, is equivalent to the following Bregman symmetric ADMM for solving \eqref{admm-sum-2}:
\begin{subequations}\label{Bregman-sym-ADMM}
\begin{align}
u^k & := \argmin_{u} \ f(u) + \frac{1}{\gamma_k}h^*(\nabla h(w^k) + \gamma_k(Mu+Nv^{k-1}-b)) \\
w^{k+\half} & := \nabla h^*(\nabla h(w^k)+\gamma_k(Mu^k+Nv^{k-1}-b)) \\
v^k & := \argmin_{v} \ g(v) + \frac{1}{\gamma_k}h^*(\nabla h(w^{k+\half}) + \gamma_k(Mu^k+Nv-b)) \\
w^{k+1} & := \nabla h^*(\nabla h(w^{k+\half})+\gamma_k(Mu^k+Nv^k-b)).
\end{align}
\end{subequations}
\end{theorem}

\begin{proof} 
The proof is similar to the proof of Theorem \ref{thm:BDRS-equiv-BADMM-linear-equality}. The BPRS \eqref{BPRS-alg} can be equivalently written as 
\begin{subequations}\label{BPRS-alg-equiv} 
\begin{align}
x^k & := J_{\gamma_k B}^h(z^k) \\
y^k & := J_{\gamma_k A}^h\circ \nabla h^*(2\nabla h(x^k) - \nabla h(z^k)) \\
z^{k+1} & := \nabla h^*(\nabla h(z^k) - 2\nabla h(x^k) + 2\nabla h(y^k)),
\end{align}
\end{subequations}
which is further equivalent to
\begin{subequations}\label{BPRS-alg-equiv-equiv} 
\begin{align}
x^k & := \argmin_x \ g^*(-N^\top x) + b^\top x + \frac{1}{\gamma_k}D_h(x,z^k) \label{BPRS-alg-equiv-equiv-1} \\
y^k & := \argmin_y \ f^*(-M^\top y) + \frac{1}{\gamma_k}D_h(y,\nabla h^*(2\nabla h(x^k) - \nabla h(z^k))) \label{BPRS-alg-equiv-equiv-2}\\
z^{k+1} & := \nabla h^*(\nabla h(z^k) - 2\nabla h(x^k) + 2\nabla h(y^k)).
\end{align}
\end{subequations}
By associating $v$ as the dual variable of \eqref{BPRS-alg-equiv-equiv-1}, and $u$ as the dual variable of \eqref{BPRS-alg-equiv-equiv-2}, similar to the proof of Theorem \ref{thm:BDRS-equiv-BADMM-linear-equality}, it can be shown that \eqref{BPRS-alg-equiv-equiv} is equivalent to 
\begin{subequations}\label{BPRS-admm-sum-2-dual-OptCond-rewrite-uv}
\begin{align}
v^k & := \argmin_v h^*(\nabla h(z^k)+\gamma_k(Nv-b)) +\gamma_k g(v) \\
x^k & := \nabla h^*(\gamma_k(Nv^k-b)+\nabla h(z^k)) \label{BPRS-admm-sum-2-dual-OptCond-rewrite-uv-x}\\
u^k & := \argmin_u h^*(2\nabla h(x^k)-\nabla h(z^k)+\gamma_k Mu) +\gamma_k f(u) \\
y^k & := \nabla h^*(2\nabla h(x^k)-\nabla h(z^k)+\gamma_k Mu^k) \label{BPRS-admm-sum-2-dual-OptCond-rewrite-uv-y}\\
z^{k+1} & := \nabla h^*(\nabla h(z^k) - 2\nabla h(x^k) + 2\nabla h(y^k)).\label{BDRS-admm-sum-2-dual-OptCond-rewrite-uv-z}
\end{align}
\end{subequations}
Reordering the updates in \eqref{BPRS-admm-sum-2-dual-OptCond-rewrite-uv}, we know that \eqref{BPRS-admm-sum-2-dual-OptCond-rewrite-uv} is equivalent to 
\begin{subequations}\label{BPRS-admm-sum-2-dual-OptCond-rewrite-uv-reorder}
\begin{align}
u^k & := \argmin_u h^*(2\nabla h(x^k)-\nabla h(z^k)+\gamma_k Mu) +\gamma_k f(u) \label{BPRS-admm-sum-2-dual-OptCond-rewrite-uv-reorder-1}\\
y^k & := \nabla h^*(2\nabla h(x^k)-\nabla h(z^k)+\gamma_k Mu^k) \label{BPRS-admm-sum-2-dual-OptCond-rewrite-uv-reorder-2}\\
z^{k+1} & := \nabla h^*(\nabla h(z^k) - 2\nabla h(x^k) + 2\nabla h(y^k)) \label{BPRS-admm-sum-2-dual-OptCond-rewrite-uv-reorder-3}\\
v^{k+1} & := \argmin_v h^*(\nabla h(z^{k+1})+\gamma_k(Nv-b)) +\gamma_k g(v) \label{BPRS-admm-sum-2-dual-OptCond-rewrite-uv-reorder-4}\\
x^{k+1} & := \nabla h^*(\gamma_k(Nv^{k+1}-b)+\nabla h(z^{k+1})). \label{BPRS-admm-sum-2-dual-OptCond-rewrite-uv-reorder-5}
\end{align}
\end{subequations}
Note that \eqref{BPRS-admm-sum-2-dual-OptCond-rewrite-uv-reorder-5} implies 
\be\label{BPRS-relation-1}
\nabla h(z^{k+1}) = \nabla h(x^{k+1}) - \gamma_k (Nv^{k+1}-b),
\ee
which immediately gives 
\be\label{BPRS-relation-2}
\nabla h(z^{k}) = \nabla h(x^{k}) - \gamma_k (Nv^{k}-b).
\ee
Combining \eqref{BPRS-admm-sum-2-dual-OptCond-rewrite-uv-reorder-2} and \eqref{BPRS-admm-sum-2-dual-OptCond-rewrite-uv-reorder-3} yields 
\be\label{BPRS-relation-3}
\nabla h(z^{k+1}) = \nabla h(y^k) + \gamma_k Mu^k,
\ee
which together with \eqref{BPRS-relation-1} gives 
\be\label{BPRS-relation-4}
\nabla h(x^{k+1}) = \nabla h(y^k) + \gamma_k (Mu^k+Nv^{k+1}-b).
\ee
Moreover, \eqref{BPRS-admm-sum-2-dual-OptCond-rewrite-uv-reorder-2} implies
\be\label{BPRS-relation-5}
\nabla h(y^{k}) = 2\nabla h(x^k) - \nabla h(z^k) + \gamma_k Mu^k = \nabla h(x^k) + \gamma_k (Mu^k+Nv^{k}-b),
\ee
where the second equality is due to \eqref{BPRS-relation-2}.
Substituting \eqref{BPRS-relation-3} into \eqref{BPRS-admm-sum-2-dual-OptCond-rewrite-uv-reorder-4}, \eqref{BPRS-relation-2} into \eqref{BPRS-admm-sum-2-dual-OptCond-rewrite-uv-reorder-1}, and combining with \eqref{BPRS-relation-5} and \eqref{BPRS-relation-4}, we know that \eqref{BPRS-admm-sum-2-dual-OptCond-rewrite-uv-reorder} is equivalent to 

\begin{align*}
u^k & := \argmin_u h^*(\nabla h(x^k)+\gamma_k (Mu+Nv^k-b)) +\gamma_k f(u) \\
y^k & := \nabla h^*(\nabla h(x^k)+\gamma_k (Mu^k+Nv^k-b)) \\ 
v^{k+1} & := \argmin_v h^*(\nabla h(y^{k})+\gamma_k(Mu^k+Nv-b)) +\gamma_k g(v) \\
x^{k+1} & := \nabla h^*(\nabla h(y^k) + \gamma_k(Mu^k+Nv^{k+1}-b)),
\end{align*}

which is exactly the same as the Bregman symmetric ADMM \eqref{Bregman-sym-ADMM}. 
\end{proof}

For the linearly inequality constrained problem \eqref{admm-sum-2-inequality-constraint}, the BPRS with $h$ being the Boltzmann-Shannon entropy leads to a symmetric version of ADEMM. 

\begin{theorem}
For $\gamma_k>0$, the BPRS \eqref{BPRS-alg} for solving the dual of \eqref{admm-sum-2-inequality-constraint} with $h$ being the Boltzmann-Shannon entropy and $A(x) = -M\partial f^*(-M^\top x)$ and $B(x) = -N \partial g^*(-N^\top x) + b$, is equivalent to the following Bregman symmetric ADEMM for solving \eqref{admm-sum-2-inequality-constraint}:
\begin{subequations}\label{sym-ADEMM-alg}
\begin{align}
u^k & := \argmin_{u} f(u) + \frac{1}{\gamma_k}\sum_{j=1}^m w_j^k \psi(\gamma_k(M_j^\top u+N_j^\top v^{k-1}-b_j)) \label{sym-ADEMM-alg-1} \\ 
w_j^{k+\half} & := w_j^k e^{\gamma_k(M_j^\top u^k+N_j^\top v^k-b_j)}, \ j=1,\ldots,m.\label{sym-ADEMM-alg-2}\\
v^k & := \argmin_{v} g(v) + \frac{1}{\gamma_j^k}\sum_{j=1}^m w_j^{k+\half}\psi(\gamma_k(M_j^\top u^k+N_j^\top v-b_j)) \label{sym-ADEMM-alg-3} \\ 
w_j^{k+1} & := w_j^{k+\half} e^{\gamma_k(M_j^\top u^k+N_j^\top v^k-b_j)}, \ j=1,\ldots,m,\label{sym-ADEMM-alg-4}
\end{align}
\end{subequations}
where $\psi$ is defined in \eqref{def-psi}.
\end{theorem}

\begin{proof}
The proof is very similar to the proof of Theorem \ref{thm:BPRS-equiv-symBADMM}. We thus omit it for brevity. 
\end{proof}

Another closely related splitting method is the double-backward (also called backward-backward) method. The double-backward method was proposed by Passty \cite{Passty-1979} and later studied by many others including Combettes \cite{Combettes-2004}. The Bregman double-backward method for solving \eqref{sum-2-operator} is given by
\be\label{BDB-alg}
z^{k+1} := J_{\gamma_kA}^hJ_{\gamma_kB}^h(z^k).
\ee
The theoretical analysis of BDBM has been mainly studied for solving the feasibility problem, i.e., when $A$ and $B$ are both normal cone operators. In this case, the monotone inclusion problem reduces to 
\be\label{prob-feasibility}
\mbox{ Find } x, \ \st, \ x \in \X_1 \bigcap \X_2,
\ee
where $\X_1$ and $\X_2$ are convex sets, and $A$ and $B$ are normal cone operators of $\X_1$ and $\X_2$, respectively. Under the assumption that $A^{-1}(0)\cap B^{-1}(0)$ is not empty, the weak convergence of BDBM \eqref{BDB-alg} was proved by Reich \cite{Reich-weak-convergence-alternating-Bregman}. In fact, note that both $J_{\gamma_kB}^h$ and $J_{\gamma_kA}^h$ are quasi-Bregman nonexpansive (QBNE, see \cite{Sabach-products-resolvents-2011}), and thus their product is also QBNE \cite{Sabach-products-resolvents-2011}. Therefore, we have $\forall u \in A^{-1}(0) \bigcap B^{-1}(0) = \Fix(J_{\gamma_kA}^h) \bigcap \Fix(J_{\gamma_kB}^h)$,
\be\label{BDBM-decrease}
D_h(u,x_{k+1}) = D_h(u,J_{\gamma_kB}^h J_{\gamma_kA}^h(x_k)) \leq D_h(u,x_k) \leq D_h(u,x_0),
\ee
where $\Fix(A)$ denotes the fixed point set of $A$. Therefore, \eqref{BDBM-decrease} implies that $D_h(u,x_k)$ is convergent.

For general $A$ and $B$, however, it may not be reasonable to assume that $A^{-1}(0)\cap B^{-1}(0)$ is not empty, and thus the convergence result in \cite{Reich-weak-convergence-alternating-Bregman} does not apply. This is indeed the case as we will see later for the optimal transport problem in Section \ref{sec:OT}. 

\section{Applications to Discrete Optimal Transport}\label{sec:OT}

Optimal transport has found many important applications in machine learning and data science recently \cite{villani2008optimal,peyre2019computational}. In the case of discrete probability measures, one is given two sets of finite number atoms, $\{y_1, y_2, \ldots, y_n\} \subset \br^d$ and  $\{z_1, z_2, \ldots, z_n\} \subset \br^d$, and two probability distributions $\mu_n = \sum_{i=1}^n r_i\delta_{y_i}$ and $\nu_n = \sum_{j=1}^n c_j\delta_{z_j}$. Here $r = (r_1, r_2, \ldots, r_n)^\top \in \Delta^n$ and $c = (c_1, c_2, \ldots, c_n)^\top \in \Delta^n$,
$\Delta^n$ denotes the probability simplex in $\br^n$ and $\delta_y$ denotes the Dirac delta function at $y$. The optimal transport between $\mu_n$ and $\nu_n$ is obtained by solving the following problem:
\be \label{OT}
 \min_{X} \langle C, X \rangle,  \ \st, \ X\onebf = r, X^\top\onebf = c, X\geq 0,
\ee
where $\onebf$ denotes the $n$-dimensional all-one vector, $C \in\br^{n\times n}$ is the cost matrix whose $(i,j)$-th component is $C_{ij} =  \|y_i - z_j\|^2$. Note that \eqref{OT} is a linear program (LP) and can be solved by off-the-shelf LP solvers. However, \eqref{OT} appearing in real applications can be very large, and classical LP solvers may suffer scalability issues. In \cite{cuturi2013sinkhorn}, Cuturi suggested to adopt the algorithm proposed by Sinkhorn and Knopp \cite{Sinkhorn-Knopp-1967} to solve the following approximation of \eqref{OT}:
\be \label{OT-reg}
 \min_{X} \langle C, X \rangle + \eta h(X), \ \st, \ X\onebf = r, X^\top\onebf = c, 
\ee
where $\eta>0$ is a penalty parameter, $h(X)$ is the Boltzmann-Shannon entropy, and for matrix $X$ it is defined as:
$h(X) = \sum_{ij} X_{ij}(\log X_{ij}-1)$. That is, an entropy penalty term is added to the objective function with a penalty parameter $\eta$. The advantage of the penalized problem is that the nonnegativity constraint $X\geq 0$ is no longer needed because it is implicitly enforced by the entropy function $h(X)$. The algorithm proposed in \cite{Sinkhorn-Knopp-1967} (from now on, we call it Sinkhorn's algorithm) solves the dual problem \eqref{OT-reg} using a block minimization algorithm. More specifically, using $\alpha$ and $\beta$ to denote the Lagrange multipliers associated with the two linear equality constraints of \eqref{OT-reg}, the dual problem of \eqref{OT-reg} can be written as:
\begin{subequations}\label{OT-reg-dual}
\begin{align}
& \max_{\alpha,\beta}\min_X \ \langle C, X \rangle + \eta h(X) - \langle \alpha, X\onebf - r\rangle - \langle \beta, X^\top\onebf - c \rangle \label{OT-reg-dual-1}\\
= &\max_{\alpha,\beta} \ \langle \alpha,r\rangle + \langle \beta,c \rangle -\eta\sum_{ij}e^{\frac{1}{\eta}(\alpha_i+\beta_j-C_{ij})}.\label{OT-reg-dual-2}
\end{align}
\end{subequations}
Note that the $X$-minimization problem in \eqref{OT-reg-dual-1} has a closed-form optimal solution given by
\[X_{ij} = e^{\frac{1}{\eta}(\alpha_i+\beta_j-C_{ij})}, \quad i,j = 1,\ldots,n.\]
By letting $K_{ij} = e^{-C_{ij}/\eta}$, $u_i = e^{\alpha_i/\eta}$, and $v_j = e^{\beta_j/\eta}$, \eqref{OT-reg-dual-2} can be equivalently written as:
\be\label{sinkhorn-uv-problem}
\min_{u,v} \sum_{ij} u_i K_{ij} v_j - \sum_i r_i \log u_i - \sum_j c_j \log v_j.
\ee
The Sinkhorn's algorithm \cite{Sinkhorn-Knopp-1967,cuturi2013sinkhorn} solves \eqref{sinkhorn-uv-problem} by alternatingly minimizing $u$ and $v$ with the other variable being fixed, i.e., 
\begin{subequations}\label{sinkhorn-alg}
\begin{align}
u^k & := \argmin_u \ \sum_{ij} u_i K_{ij} v_j^{k-1} - \sum_i r_i \log u_i  \\
v^k & := \argmin_v \ \sum_{ij} u_i^k K_{ij} v_j - \sum_j c_j \log v_j.
\end{align}
\end{subequations}
It turns out that both subproblems in \eqref{sinkhorn-alg} admit closed-form solutions, and the Sinkhorn's algorithm can be written more compactly as:
\begin{subequations}\label{sinkhorn-alg-simpler}
\begin{align}
u^k & := r./(Kv^{k-1})   \\
v^k & := c./(K^\top u^k) .
\end{align}
\end{subequations}
The Sinkhorn's algorithm can be implemented very efficiently because the computations in \eqref{sinkhorn-alg-simpler} are simple. One drawback of the Sinkhorn's algorithm is that it is very difficult to tune the parameter $\eta$. Ideally, one wants a small $\eta$ so that the penalized problem \eqref{OT-reg} is close to the original problem \eqref{OT}. However, small $\eta$ will cause numerical instability of the Sinkhorn's algorithm. It is usually not clear how small $\eta$ can be without causing numerical issues. Moreover, the Sinkhorn's algorithm only solves the regularized problem \eqref{OT-reg}, not the original OT problem \eqref{OT}. 

Now we discuss how to use our BDRS (or ADEMM) to solve the original OT \eqref{OT}. Note that \eqref{OT} can be written in the form of \eqref{sum-2-operator} by defining 
\be\label{def-OT-AB}
A(X) = C + \partial \bone(X\onebf=r) + \partial\bone(X\geq 0), \quad B(X) = \partial \bone(X^\top \onebf=c) + \partial\bone(X\geq 0).
\ee
Now the BDRS \eqref{BDRS-alg} can be written as (for the ease of comparison with Sinkhorn's algorithm, we choose $\gamma_k$ to be a constant $1/\eta$):
\begin{subequations}\label{BDRS-alg-OT}
\begin{align}
X^k & := \argmin_X \ D_h(X,Z^k), \ \st, \ X^\top \onebf=c \label{BDRS-alg-OT-1}\\ 
Y^k & := \argmin_Y \ \langle C,Y\rangle + \eta D_h(Y, X^k\circ X^k./Z^k), \ \st, \ Y\onebf=r \label{BDRS-alg-OT-2}\\ 
Z^{k+1} & := Z^k \circ Y^k./ X^k.\label{BDRS-alg-OT-3}
\end{align}
\end{subequations}
The dual problem of the OT \eqref{OT} is given by:
\be\label{OT-dual}
\min_{\alpha, \beta} \ -r^\top \alpha - c^\top \beta, \ \st, \ \alpha_i + \beta_j \leq C_{ij}, i, j = 1,\ldots,n.
\ee
The ADEMM \eqref{ADEMM-alg} for solving \eqref{OT-dual} is given by:
\begin{subequations}\label{ADEMM-alg-OT}
\begin{align}
\alpha^k & :=  \argmin_\alpha \ -r^\top \alpha + \eta\sum_{ij}X_{ij}^k e^{\frac{1}{\eta}(\alpha_i+\beta_j^{k-1}-C_{ij})} \\
\beta^k & := \argmin_\beta \ - c^\top \beta + \eta\sum_{ij}X_{ij}^k e^{\frac{1}{\eta}(\alpha_i^k+\beta_j-C_{ij})} \\
X_{ij}^{k+1} & := X_{ij}^ke^{\frac{1}{\eta}(\alpha_i^k+\beta_j^k-C_{ij})}, \ i,j=1,\ldots,n.
\end{align}
\end{subequations}
By denoting $u = e^{\alpha/\eta}$, $v= e^{\beta/\eta}$, and $K_{ij} = e^{-C_{ij}/\eta}$, \eqref{ADEMM-alg-OT} can be further rewritten as 
\begin{subequations}\label{ADEMM-alg-OT-rewrite-uv}
\begin{align}
u^k & := \argmin_u \ \sum_{ij}u_i X_{ij}^kK_{ij}v_j^{k-1} -\sum_i r_i\log u_i \label{ADEMM-alg-OT-rewrite-uv-1}\\
v^k & := \argmin_v \ \sum_{ij}u_i^k X_{ij}^kK_{ij}v_j - \sum_j c_j\log v_j \label{ADEMM-alg-OT-rewrite-uv-2}\\
X_{ij}^{k+1} & := u_i^kX_{ij}^kK_{ij}v_j^k, \ i,j=1,\ldots,n.
\end{align}
\end{subequations}
Similar to \eqref{sinkhorn-alg}, the two subproblems \eqref{ADEMM-alg-OT-rewrite-uv-1} and \eqref{ADEMM-alg-OT-rewrite-uv-2} admit closed-form solutions, and \eqref{ADEMM-alg-OT-rewrite-uv} can be equivalently written as:
\begin{subequations}\label{ADEMM-alg-OT-rewrite-uv-simpler}
\begin{align}
u^k & := r./((X^k\circ K)v^{k-1}) \label{ADEMM-alg-OT-rewrite-uv-simpler-1}\\
v^k & := c./((X^k\circ K)^\top u^k) \label{ADEMM-alg-OT-rewrite-uv-simpler-2}\\
X_{ij}^{k+1} & := u_i^kX_{ij}^kK_{ij}v_j^k, \ i,j=1,\ldots,n.\label{ADEMM-alg-OT-rewrite-uv-simpler-3}
\end{align}
\end{subequations}
Now comparing ADEMM \eqref{ADEMM-alg-OT-rewrite-uv-simpler} with the Sinkhorn's algorithm \eqref{sinkhorn-alg-simpler}, we see when updating $u^k$ and $v^k$, \eqref{ADEMM-alg-OT-rewrite-uv-simpler} replaces matrix $K$ in \eqref{sinkhorn-alg-simpler} by matrix $X^k\circ K$, which is no longer a constant matrix. The matrix $X^{k+1}$ is then updated by \eqref{ADEMM-alg-OT-rewrite-uv-simpler-3}. 

\begin{remark}
Our BDRS \eqref{ADEMM-alg-OT-rewrite-uv-simpler} for solving the OT probelm can be a much better algorithm than the Sinkhorn's algorithm, because we do not require $\eta$ to be close to 0, and thus resolve the issue of numerical instability of the Sinkhorn's algoirthm. We also see that the computations in \eqref{ADEMM-alg-OT-rewrite-uv-simpler} are very simple and can be done in parallel as illustrated in \cite{cuturi2013sinkhorn} for the Sinkhorn's algorithm. Thus our BDRS can be a much better algorithm than the classical ADMM for solving \eqref{OT} which requires projections onto the probability simplex in each iteration. 
\end{remark}

We now discuss the BDBM \eqref{BDB-alg} for solving the OT problem \eqref{OT}.
We note that $\Fix(J_{\gamma A}^h J_{\gamma B}^h) \not\subset (A+B)^{-1}(0)$, although $\Fix(R_{\gamma A}^h R_{\gamma B}^h) \subset (A+B)^{-1}(0)$. Therefore, even though we can find a fixed point of $J_{\gamma A}^h J_{\gamma B}^h$ using the BDBM with a constant $\gamma_k=\gamma$, the solution we find may not solve the OT problem \eqref{OT}. However, we have the following result regarding the fixed point of $J_{\gamma_k A}^h J_{\gamma_k B}^h$. 

\begin{theorem}\label{thm:BDBM-OT}
For the OT problem \eqref{OT} with $A$ and $B$ defined in \eqref{def-OT-AB} and $h$ being the Boltzmann-Shannon entropy, it holds that 
\[\lim_{\gamma_k\to 0} \Fix(J_{\gamma_k A}^h J_{\gamma_k B}^h) \subset (A+B)^{-1}(0).\]
\end{theorem}

\begin{proof}
Assume $x\in \Fix(J_{\gamma_k A}^h J_{\gamma_k B}^h)$, that is 
\[ x = J_{\gamma_k A}^h J_{\gamma_k B}^h x,\]
which can be equivalently written as 
\[x = (\nabla h + \gamma_k A)^{-1}\circ \nabla h \circ (\nabla h + \gamma_k B)^{-1}\circ \nabla h(x).\]
This is further equivalent to (note that $\nabla h(x) = \log x$ and $\nabla h^*(x) = e^x$):
\begin{align*}
& \quad (\nabla h + \gamma_k B)\nabla h^*(\nabla h + \gamma_k A)x = \nabla h(x) \\
\Longleftrightarrow & \quad (I + \gamma_k B\nabla h^*) (\nabla h (x) + \gamma_k A(x)) = \nabla h(x) \\
\Longleftrightarrow & \quad (I + \gamma_k B\nabla h^*) \log (x^k\cdot e^{\gamma_k A(x)}) = \nabla h(x) \\
\Longleftrightarrow & \quad \log (x\cdot e^{\gamma_k A(x)}) + \gamma_k B(x\cdot e^{\gamma_k A(x)}) = \nabla h(x) \\
\Longleftrightarrow & \quad \log (x\cdot e^{\gamma_k A(x)} \cdot e^{\gamma_k B(x\cdot e^{\gamma_k A(x)})}) = \log(x) \\
\Longleftrightarrow & \quad \gamma_k A(x) + \gamma_k B(x\cdot e^{\gamma_k A(x)}) = 0 \\
\Longleftrightarrow & \quad A(x) + B(x\cdot e^{\gamma_k A(x)}) = 0.
\end{align*}
By letting $\gamma_k \rightarrow 0$, we have $A(x)+ B(x) = 0$, i.e., $x\in (A+B)^{-1}(0)$.
\end{proof}

Theorem \ref{thm:BDBM-OT} shows that if we let $\gamma_k\rightarrow 0$, then the BDBM \eqref{BDB-alg} solves the OT problem \eqref{OT}.

\section{Concluding Remarks}\label{sec:conclusion}

In this paper, we proposed several new algorithms for solving monotone operator inclusion problem and convex minimization problems. The algorithms include BDRS, BPRS, Bregman ADMM, and ADEMM. We discussed their connections with existing algorithms in the literature. We also discussed how to apply our algorithms to solve the discrete optimal transport problem. We proved the convergence of the algorithms under certain assumptions, though we point out that one assumption does not apply to the OT problem. We leave it as a future work to prove the convergence of the proposed algorithms under the most general setting. 

\bibliographystyle{plain}
\bibliography{mach_learn,BDRS-and-manifold}

\appendix 
\section{The Convergence Analysis of BPRS}\label{appendix:BPRS}

Note that an interesting observation to the operators defined in Definition \ref{def-Bregman-operators} is that 
\[R_T^h = F_T^h\circ J_T^h.\]

Therefore the BPRS \eqref{BPRS-alg} is equivalent to
\[
z^{k+1} = F_{\gamma_kA}^h\circ J_{\gamma_kA}^h \circ F_{\gamma_kB}^h\circ J_{\gamma_kB}^h (z^k).
\]
Let $x^k =  J_{\gamma_kA}^h \circ F_{\gamma_kB}^h\circ J_{\gamma_kB}^h (z^k)$, then the above equation becomes 
\be\label{BPRS-x-special}
x^{k+1} = J_{\gamma_kA}^h\circ F_{\gamma_kB}^h \circ J_{\gamma_kB}^h\circ F_{\gamma_kA}^h (x^k).
\ee
Therefore, the BPRS is the composition of two Bregman forward-backward operators. Moreover it can be shown that 
\be\label{BPRS-special-fixed-set}
\Fix(J_{\gamma_kA}^h\circ F_{\gamma_kB}^h \circ J_{\gamma_kB}^h\circ F_{\gamma_kA}^h) = (A+B)^{-1}(0).
\ee

\subsection{The convergence of BPRS when both $f$ and $g$ are relatively smooth.}

In this section, we provide a convergence analysis of BPRS when both 
$f$ and $g$ are relative smooth functions with respect to $D_h$, i.e.,
\begin{align*}
f(x) \leq f(y) + \langle \nabla f(y), x-y \rangle + LD_h(x,y) \\
g(x) \leq g(y) + \langle \nabla g(y), x-y \rangle + LD_h(x,y), 
\end{align*}
where $L>0$ is the relative smoothness parameter. 
We consider the smooth problem 
\be\label{proof-problem-special}
\min_x \ F(x) = f(x) + g(x) \Longleftrightarrow  0 \in  A(x) + B(x),
\ee
where $A = \nabla f$ and $B = \nabla g$.
Our convergence result is summarized in Theorem \ref{thm:BPRS-convergence-smooth}.

\begin{theorem}\label{thm:BPRS-convergence-smooth}
    Assume $f$ and $g$ are relative smooth functions with respect to $D_h$ with parameter $L$. BPRS with $\gamma_k=\gamma\leq 1/L$ for solving \eqref{proof-problem-special} globally converges to $(A+B)^{-1}(0)$.
\end{theorem}
\begin{proof}
BPRS with $\gamma_k=\gamma\leq 1/L$ for solving \eqref{proof-problem-special} can be rewritten as
\begin{subequations}\label{BPRS-alg-special}
\begin{align}
y^k & := \argmin_y \ f(x^k) + \langle \nabla f(x^k), y-x^k\rangle + g(y) + \frac{1}{\gamma}D_h(y,x^k) \\
x^{k+1} & := \argmin_x \ g(y^k) + \langle \nabla g(y^k), x-y^k\rangle + f(x) + \frac{1}{\gamma}D_h(x,y^k).
\end{align}
\end{subequations}
It is easy to obtain that:
\begin{align*}
F(y^k) \leq & F(x^k) - \frac{1}{\gamma} D_h(x^k,y^k) \\
F(x^{k+1}) \leq & F(y^k) - \frac{1}{\gamma} D_h(y^k,x^{k+1}),
\end{align*}
which further yields
\begin{align*} 
F(x^{k+1}) \leq F(x^k) - \frac{1}{\gamma} D_h(x^k,y^k) - \frac{1}{\gamma} D_h(y^k,x^{k+1}).
\end{align*}
This inequality shows that $F(x^k)$ monotonically decreases. Moreover, by telescoping sum, we obtain:
\[\frac{1}{\gamma}\sum_{k=0}^N \left(D_h(x^k,y^k) + D_h(y^k,x^{k+1}) \right)\leq F(x^0) -F(x^*).\]
This indicates 
\[\lim_{k\to \infty} D_h(x^k,y^k) = 0, \mbox{ and } \lim_{k\to \infty} D_h(y^k,x^{k+1}) = 0.\]
This further implies
\[\lim_{k\to \infty} \|x^k-y^k\|=0, \quad \lim_{k\to \infty} \|x^{k+1}-y^k\|=0, \quad \lim_{k\to \infty} \|x^k-x^{k+1}\|=0.\]
This shows that $\{x^k\}$ is a Cauchy sequence. Therefore, $\{x^k\}$ is bounded  and convergent. From \eqref{BPRS-x-special} and \eqref{BPRS-special-fixed-set} we know that  $\{x^k\}$ converges to $(A+B)^{-1}(0)$.
\end{proof}

\subsection{The convergence of BPRS when both $f$ and $g$ are non-smooth functions.}\label{sec:appendix-BPRS-nonsmooth}

Note that even when both $f$ and $g$ are nonsmooth, we can still apply \eqref{BPRS-x-special} with $A = \partial f$ and $B = \partial g$. In this case, \eqref{BPRS-x-special} reduces to {(with $x^0:=J_{\gamma_kA}^h \circ F_{\gamma_kB}^h\circ J_{\gamma_kB}^h (z^0)$)}
\be\label{BPRS-alternating-prox-subgrad}
x^{k+1} = (\nabla h+\gamma_k \partial f)^{-1} (\nabla h - \gamma_k \partial g) (\nabla h+\gamma_k \partial g)^{-1} (\nabla h - \gamma_k \partial f) (x^k).
\ee

In this section we prove the convergence of the BPRS for solving 
\be\label{prove-BPRS-prob-sum-1}\min_x \ f(x) + g(x)\ee
when both $f$ and $g$ are nonsmooth functions, under the assumption that 
$\textrm{im}(\nabla h^*) \subset \dom f \cap \dom g$. This assumption was used in \cite{Bot-Bregman-prox-subgrad}. 
Note that according to \eqref{BPRS-alternating-prox-subgrad}, the BPRS for solving \eqref{prove-BPRS-prob-sum-1} can be written as
\begin{subequations}\label{BPRS-alt-prox-subgrad-alg}
\begin{align}
\bar{x}^k & := \nabla h^*(\nabla h(w^k) - \gamma_k\partial f(w^k)) \label{BPRS-alt-prox-subgrad-alg-1}\\
x^{k+1}  & := \argmin_x \ g(x) + \frac{1}{\gamma_k}D_h(x,\bar{x}^k) \label{BPRS-alt-prox-subgrad-alg-2}\\
\bar{w}^{k+1} &:= \nabla h^*(\nabla h(x^{k+1}) - \gamma_k\partial g(x^{k+1})) \label{BPRS-alt-prox-subgrad-alg-3}\\
w^{k+1} & := \argmin_w \ f(w) + \frac{1}{\gamma_k}D_h(w,\bar{w}^{k+1}).\label{BPRS-alt-prox-subgrad-alg-4}
\end{align}
\end{subequations}
Note that this is an alternating Bregman proximal subgradient method. Bot and Bohm \cite{Bot-Bregman-prox-subgrad} analyzed the convergence of Bregman proximal subgradient method, which consists only one step of \eqref{BPRS-alt-prox-subgrad-alg}. 

We will use the {\emph{Three point identity:}}
\be\label{three-point} 
D_h(x,y) + D_h(y,z) = D_h(x,z) - \langle\nabla h(y)-\nabla h(z),x-y \rangle.
\ee

Our main result of the convergence of BPRS is summarized in Theorem \ref{thm:BPRS-convergence-nonsmooth}.

\begin{theorem}\label{thm:BPRS-convergence-nonsmooth}
    Assume $\textrm{im}(\nabla h^*) \subset \dom f \cap \dom g$, $\|\partial f\|_\infty\leq G$, $\|\partial g\|_\infty\leq G$, $h$ is $\sigma$-strongly convex over the simplex, and 
    $\nabla h^*$ is Lipschitz continuous. Then BPRS \eqref{BPRS-alt-prox-subgrad-alg} with $\gamma_k=1/\sqrt{k}$ for solving \eqref{prove-BPRS-prob-sum-1} converges sublinearly. 
\end{theorem}

\begin{proof} 
First, we have
\begin{align}
D_h(y,\bar{x}^k) & \overset{\eqref{three-point}}{=\joinrel=} D_h(y,w^k) - D_h(\bar{x}^k,w^k) - \langle \nabla h(\bar{x}^k) - \nabla h(w^k), y-\bar{x}^k\rangle \nonumber\\
& = D_h(y,w^k) - D_h(\bar{x}^k,w^k) + \gamma_k \langle\partial f(w^k), y-\bar{x}^k\rangle \nonumber\\
& =  D_h(y,w^k) - D_h(\bar{x}^k,w^k) + \gamma_k \langle \partial f(w^k), y-w^k\rangle - \gamma_k \langle \partial f(w^k), \bar{x}^k-w^k\rangle \nonumber\\
& \leq D_h(y,w^k) - D_h(\bar{x}^k,w^k) + \gamma_k (f(y) - f(w^k)) + \gamma_k \|\partial f(w^k)\|_\infty\|\bar{x}^k-w^k\|_1\nonumber \\
& \leq D_h(y,w^k) - D_h(\bar{x}^k,w^k) + \gamma_k (f(y) - f(w^k)) + \frac{\gamma_k^2}{\sigma}\|\partial f(w^k)\|_\infty^2 + \frac{\sigma}{4}\|\bar{x}^k-w^k\|_1^2 \nonumber\\
& \overset{(*)}{\leq} D_h(y,w^k) - D_h(\bar{x}^k,w^k) + \gamma_k (f(y) - f(w^k)) + \frac{\gamma_k^2}{\sigma}\|\partial f(w^k)\|_\infty^2 + \frac{1}{2}D_h(\bar{x}^k,w^k) \nonumber\\
& \overset{(**)}{\leq} D_h(y,w^k) - \frac{1}{2}D_h(\bar{x}^k,w^k) + \gamma_k (f(y) - f(w^k)) + \frac{\gamma_k^2}{\sigma}G^2 \label{proof-bound-1},
\end{align} 
where (*) is due to the $\sigma$-strong convexity of $h$ over the simplex (also note that both $\bar{x}^k$ and $w^k$ are on simplex), in (**) we used $\|\partial f\|_\infty\leq G$. 

Now, the optimality condition of \eqref{BPRS-alt-prox-subgrad-alg-2} is:
\[0\in\gamma_k\partial g(x^{k+1}) + \nabla h(x^{k+1}) - \nabla h(\bar{x}^k),\]
which implies
\[\gamma_k(g(y) - g(x^{k+1})) \geq \langle \nabla h(\bar{x}^k)-\nabla h(x^{k+1}), y - x^{k+1}\rangle.\]
Using the three-point identity we have
\[\gamma_k(g(y) - g(x^{k+1})) \geq D_h(y,x^{k+1}) + D_h(x^{k+1},\bar{x}^k) - D_h(y,\bar{x}^k),\]
or, equivalently, 
\be\label{proof-bound-2}\gamma_k(g(x^{k+1})-g(y)) + D_h(y,x^{k+1}) \leq D_h(y,\bar{x}^k)- D_h(x^{k+1},\bar{x}^k).
\ee
Combining \eqref{proof-bound-1} and \eqref{proof-bound-2}, we have
\[
\gamma_k(g(x^{k+1})-g(y)) + \gamma_k (f(w^k) - f(y)) + D_h(y,x^{k+1}) \leq D_h(y,w^k) + \frac{\gamma_k^2}{\sigma}G^2 - D_h(x^{k+1},\bar{x}^k) - \frac{1}{2}D_h(\bar{x}^k,w^k).
\]
By adding and subtracting $f(x^{k+1})$, we obtain 

\begin{align*}
&\gamma_k(f(x^{k+1})+g(x^{k+1})-f(y)- g(y)) + \gamma_k (f(w^k) - f(x^{k+1})) + D_h(y,x^{k+1}) \\
\leq &  D_h(y,w^k) + \frac{\gamma_k^2}{\sigma}G^2 - D_h(x^{k+1},\bar{x}^k) - \frac{1}{2}D_h(\bar{x}^k,w^k), 
\end{align*}
which immediately implies
\begin{align*}
&\gamma_k(f(x^{k+1})+g(x^{k+1})-f(y)- g(y))  + D_h(y,x^{k+1}) \\
\leq & D_h(y,w^k) + \frac{\gamma_k^2}{\sigma}G^2 - D_h(x^{k+1},\bar{x}^k) - \frac{1}{2}D_h(\bar{x}^k,w^k) + \gamma_k G \|w^k - x^{k+1}\|_1 \\
\overset{(*)} \leq &D_h(y,w^k) + \frac{\gamma_k^2}{\sigma}G^2 - D_h(x^{k+1},\bar{x}^k) - \frac{1}{2}D_h(\bar{x}^k,w^k) + \gamma_k G \|w^k - \bar{x}^k\|_1 + \gamma_kG\|\bar{x}^k-x^{k+1}\|_1 \\
\overset{(**)}\leq &D_h(y,w^k) + \frac{\gamma_k^2}{\sigma}G^2 - D_h(x^{k+1},\bar{x}^k) - \frac{1}{2}D_h(\bar{x}^k,w^k) + \gamma_k G \|w^k - \bar{x}^k\|_1 + \frac{\gamma_k^2G^2}{2\sigma} + D_h(x^{k+1},\bar{x}^k) \\
= & D_h(y,w^k) + \frac{\gamma_k^2}{\sigma}G^2 - \frac{1}{2}D_h(\bar{x}^k,w^k) + \gamma_k G \|w^k - \bar{x}^k\|_1 + \frac{\gamma_k^2G^2}{2\sigma} \\
\overset{(***)}\leq & D_h(y,w^k) + \frac{\gamma_k^2}{\sigma}G^2 - \frac{1}{2}D_h(\bar{x}^k,w^k) + \frac{\gamma_k^2G^2}{\sigma} + \frac{\sigma}{4} \|w^k - \bar{x}^k\|_1^2 + \frac{\gamma_k^2G^2}{2\sigma} \\
\overset{(****)}{\leq} & D_h(y,w^k) + \frac{\gamma_k^2}{\sigma}G^2 + \frac{\gamma_k^2G^2}{\sigma} + \frac{\gamma_k^2G^2}{2\sigma} \\
= & D_h(y,w^k) + \frac{5\gamma_k^2}{2\sigma}G^2,
\end{align*}
where (*) is the triangle inequality, (**) is due to Young's inequality and the $\sigma$-strong convexity of $h$, (***) is due to Young's inequality, and (****) is due to the $\sigma$-strong convexity of $h$. Note that so far we only dealt with \eqref{BPRS-alt-prox-subgrad-alg-1}-\eqref{BPRS-alt-prox-subgrad-alg-2}, and we obtained 
\be\label{key-inequality-1}
\gamma_k(f(x^{k+1})+g(x^{k+1})-f(y)- g(y))  + D_h(y,x^{k+1}) \leq D_h(y,w^k) + \frac{5\gamma_k^2}{2\sigma}G^2. 
\ee
Similarly, for \eqref{BPRS-alt-prox-subgrad-alg-3}-\eqref{BPRS-alt-prox-subgrad-alg-4}, we have 
\be\label{key-inequality-2}
\gamma_k(f(w^{k+1})+g(w^{k+1})-f(y)- g(y))  + D_h(y,w^{k+1}) \leq D_h(y,x^{k+1}) + \frac{5\gamma_k^2}{2\sigma}G^2. 
\ee
Now summing up \eqref{key-inequality-1} and \eqref{key-inequality-2} we have (denote $L = \frac{5G^2}{2\sigma}$)
\be\label{important-inequality-3}
\gamma_k (f(x^{k+1}) + g(x^{k+1}) + f(w^{k+1}) + g(w^{k+1}) - 2(f(y) + g(y))) + D_h(y,w^{k+1}) \leq D_h(y,w^{k}) + 2L \gamma_k^2.
\ee
Define $z^{k+1} := (x^{k+1}+w^{k+1})/2$. From the convexity of $f$ and $g$, we have
\be\label{important-inequality-4}
\gamma_k (f(z^{k+1}) + g(z^{k+1}) - (f(y) + g(y))) + \frac{1}{2}D_h(y,w^{k+1}) \leq \frac{1}{2}D_h(y,w^{k}) + L \gamma_k^2.
\ee
Now summing up \eqref{important-inequality-4} for $k = 0,1,\ldots,N-1$, we have 
\be\label{important-inequality-5}
\sum_{k=0}^{N-1}\gamma_k (f(z^{k+1}) + g(z^{k+1}) - (f(y) + g(y))) \leq \frac{1}{2}D_h(y,w^{0}) + L \sum_{k=0}^{N-1}\gamma_k^2.
\ee
Therefore, we have
\be\label{important-inequality-6}
\min_{0\leq k\leq N-1}\left(f(z^{k+1}) + g(z^{k+1})\right) - (f(y) + g(y)) \leq \frac{D_h(y,w^{0})}{2 \sum_{k=0}^{N-1}\gamma_k } + \frac{L \sum_{k=0}^{N-1}\gamma_k^2}{\sum_{k=0}^{N-1}\gamma_k}.
\ee
By choosing $y = x^*$ and using $\gamma_k = \frac{1}{\sqrt{k}}$, we obtain
\[
\lim_{N\rightarrow +\infty} \min_{0\leq k\leq N-1}\left(f(z^{k+1}) + g(z^{k+1})\right) = f(x^*) + g(x^*). 
\]
Moreover, the convergence rate is sublinear.
\end{proof}

\begin{remark}
    Note that we required that $\nabla h^*$ is Lipschitz continuous. This is satisfied if $h$ is the entropy function: $h(x) = \sum_i (x_i \log x_i - x_i)$. In this case, $\nabla h(x) = \log x$, and the conjugate function on probability simplex is given by:
\[
h^*(x) = \sup_{e^\top y =1} x^\top y - h(y). 
\]
It is easy to verify that the optimal $y$ is given by $y_i = \frac{e^{x_i}}{\sum_j e^{x_j}}$, and 
\be\label{def-h*}
h^*(x) = \log \sum_i e^{x_i} + 1. 
\ee
Note that this is $\log \sum\exp$ function, and it is known to have Lipschitz continuous gradient. And, we have $[\nabla h^*(x)]_i = \frac{e^{x_i}}{\sum_j e^{x_j}}$. 
\end{remark}

\section{Convergence of BDRS}\label{appendix:BDRS}

In this section we consider the same problem as in Section \ref{sec:appendix-BPRS-nonsmooth} under the same assumptions in Theorem \ref{thm:BPRS-convergence-nonsmooth}. 
Note that BDRS can be equivalently written as
\begin{subequations}\label{BDRS-rewrite-better}
\begin{align}
w^k & := \argmin_w \ f(w) + \frac{1}{\gamma_k}D_h(w,x^k)\label{BDRS-rewrite-1-better} \\
\bar{x}^k & := \nabla h^*(2\nabla h(w^k) - \nabla h(x^k))\label{BDRS-rewrite-2-better} \\
z^k & := \argmin_z \ g(z) + \frac{1}{\gamma_k}D_h(z,\bar{x}^k)\label{BDRS-rewrite-3-better} \\
y^k & := \nabla h^*(2\nabla h(z^k) - \nabla h(\bar{x}^k))\label{BDRS-rewrite-4-better} \\
x^{k+1} & := \nabla h^*\left(\frac{1}{2}\nabla h(x^k)+\frac{1}{2}\nabla h(y^k)\right).\label{BDRS-rewrite-5-better}
\end{align}
\end{subequations}

This is equivalent to the following when considering problem \eqref{prove-BPRS-prob-sum-1}: 

\begin{subequations}\label{BDRS-rewrite}
\begin{align}
w^k & := \argmin_w \ f(w) + \frac{1}{\gamma_k}D_h(w,x^k)\label{BDRS-rewrite-1} \\
\bar{x}^k & := \nabla h^*(\nabla h(w^k) - \gamma_k \partial f(w^k))\label{BDRS-rewrite-2} \\
z^k & := \argmin_z \ g(z) + \frac{1}{\gamma_k}D_h(z,\bar{x}^k)\label{BDRS-rewrite-3} \\
y^k & := \nabla h^*(\nabla h(z^k) - \gamma_k\partial g(z^k))\label{BDRS-rewrite-4} \\
x^{k+1} & := \nabla h^*\left(\frac{1}{2}\nabla h(x^k)+\frac{1}{2}\nabla h(y^k)\right).\label{BDRS-rewrite-5}
\end{align}
\end{subequations}

\begin{theorem}\label{thm:BDRS-convergence}
    Assume $\textrm{im}(\nabla h^*) \subset \dom f \cap \dom g$, $\|\partial f\|_\infty\leq G$, $\|\partial g\|_\infty\leq G$, $h$ is $\sigma$-strongly convex over the simplex, and 
    $\nabla h^*$ is Lipschitz continuous. Then BDRS \eqref{BDRS-rewrite} with $\gamma_k=1/\sqrt{k}$ for solving \eqref{prove-BPRS-prob-sum-1} converges sublinearly. 
\end{theorem}

Note that \eqref{BDRS-rewrite-1} is one step of Bregman PPA, which yields 
\be\label{BPPA-inequality}
\gamma_k (f(y) - f(x^{k+1})) \geq D_h(y,x^{k+1}) - D_h(y,x^k) + D_h(x^{k+1},x^k), \quad \forall y.
\ee
\eqref{BDRS-rewrite-2} is one step of the mirror subgradient method, which yields
\be\label{mirror-subgrad-inequality}
\gamma_k (f(x^k)-f(y)) \leq D_h(y,x^k) - D_h(y,x^{k+1}) - \frac{1}{2}D_h(x^{k+1},x^k)+ \frac{\gamma_k^2\|f'(x^k)\|_*^2}{2\sigma}, \quad \forall y.
\ee

Apply \eqref{mirror-subgrad-inequality} to \eqref{BDRS-rewrite-2} and \eqref{BDRS-rewrite-4}, we have (we ignored the constant coefficient of $\gamma_k^2$)
\begin{align*}
\gamma_k (f(w^k)-f(y)) & \leq D_h(y,w^k) - D_h(y,\bar{x}^{k}) - \frac{1}{2}D_h(\bar{x}^{k},w^k)+ \gamma_k^2 \\
\gamma_k (g(z^k)-g(y)) & \leq D_h(y,z^k) - D_h(y,y^{k}) - \frac{1}{2}D_h(y^{k},z^k)+ \gamma_k^2.
\end{align*}
Apply \eqref{BPPA-inequality} to \eqref{BDRS-rewrite-1} and \eqref{BDRS-rewrite-3}, we have
\begin{align*}
\gamma_k(f(w^k)-f(y)) &\leq D_h(y,x^k) - D_h(y,w^k) - D_h(w^k,x^k) \\
\gamma_k(g(z^k)-g(y)) & \leq D_h(y,\bar{x})-D_h(y,z^k) - D_h(z^k,\bar{x}^k).
\end{align*}
Combining these four inequalities, we have
\begin{align}
& 2\gamma_k(f(w^k)+g(z^k)-f(y)-g(y)) \label{combine-mirror-subgrad-BPPA} \\ \leq & D_h(y,x^k) - D_h(y,y^{k}) - \frac{1}{2}D_h(y^{k},z^k) - \frac{1}{2}D_h(\bar{x}^{k},w^k)- D_h(w^k,x^k) - D_h(z^k,\bar{x}^k)+ \gamma_k^2.\nonumber
\end{align}

Note that \eqref{BDRS-rewrite-5} is equivalent to 
\be\label{BDRS-rewrite-5-equiv}
\nabla h(x^{k+1}) - \nabla h(y^k) = \nabla h(x^k) - \nabla h(x^{k+1}). 
\ee
Using the 3-point identity \eqref{three-point} twice, we have
\begin{align*}
D_h(y,x^{k+1}) & = D_h(y,y^k) - D_h(x^{k+1},y^k)-\langle \nabla h(x^{k+1})-\nabla h(y^k),y-x^{k+1}\rangle \\
\langle \nabla h(x^k)-\nabla h(x^{k+1}),y-x^{k+1}\rangle & = D_h(y,x^{k+1}) + D_h(x^{k+1},x^k) - D_h(y,x^k). 
\end{align*}
Combining these two identities with \eqref{BDRS-rewrite-5-equiv}, we obtain 
\[
D_h(y,x^{k+1}) = D_h(y,y^k) - D_h(x^{k+1},y^k) -D_h(y,x^{k+1}) - D_h(x^{k+1},x^k) + D_h(y,x^k),
\]
which is equivalent to
\be\label{BDRS-last-step-3-pts-identity}
2D_h(y,x^{k+1}) = D_h(y,x^k) + D_h(y,y^k) - D_h(x^{k+1},y^k)-D_h(x^{k+1},x^k). 
\ee
Combining \eqref{combine-mirror-subgrad-BPPA} and \eqref{BDRS-last-step-3-pts-identity}, we have 
\begin{align}\label{combine-all-5}
& 2\gamma_k(f(w^k)+g(z^k)-f(y)-g(y)) +2D_h(y,x^{k+1}) \\ \leq & 2D_h(y,x^k) - \frac{1}{2}D_h(y^{k},z^k) - \frac{1}{2}D_h(\bar{x}^{k},w^k)- D_h(w^k,x^k) - D_h(z^k,\bar{x}^k) - D_h(x^{k+1},y^k)-D_h(x^{k+1},x^k) + \gamma_k^2.\nonumber
\end{align}
{Under the assumption that  $\textrm{im}(\nabla h^*)\subset \dom f\cap \dom g$}, we know $f(x^{k+1})$ and $g(x^{k+1})$ are finite values. By adding and subtracting $f(x^{k+1}) + g(x^{k+1})$ to \eqref{combine-all-5}, we get
\begin{align}
& 2\gamma_k(f(x^{k+1})+g(x^{k+1})-f(y)-g(y)) +2D_h(y,x^{k+1}) \label{final-bound-inequality}\\ \leq & 2D_h(y,x^k) - \frac{1}{2}D_h(y^{k},z^k) - \frac{1}{2}D_h(\bar{x}^{k},w^k)- D_h(w^k,x^k) - D_h(z^k,\bar{x}^k) - D_h(x^{k+1},y^k)-D_h(x^{k+1},x^k) \nonumber\\ & + 2\gamma_k(f(x^{k+1})- f(w^k)+g(x^{k+1})-g(z^k))+ \gamma_k^2. \nonumber
\end{align}
Now we only need to show that $2\gamma_k(f(x^{k+1})- f(w^k)+g(x^{k+1})-g(z^k))$ can be bounded by those negative terms on the right-hand-side. We have
\begin{align*}
& 2\gamma_k(f(x^{k+1})- f(w^k)+g(x^{k+1})-g(z^k)) \\ 
\leq & 2\gamma_k G (\|x^{k+1}-w^k\|_1 + \|x^{k+1}-z^k\|_1) \\
\leq & 2\gamma_k G(\|x^{k+1}-w^k\|_1 +\|z^k-w^k\|_1+\|x^{k+1}-w^k\|_1) \\
= & 2\gamma_k G(2\|x^{k+1}-w^k\|_1 +\|z^k-w^k\|_1) \\
\leq & 2\gamma_k G(2\|x^k-x^{k+1}\|_1+2\|w^k-x^k\|_1+\|z^k-w^k\|_1)\\
\leq & 2\gamma_k G(2\|x^k-x^{k+1}\|_1+2\|w^k-x^k\|_1+\|z^k-\bar{x}^k\|_1+\|\bar{x}^k-w^k\|_1)
\end{align*}
Now apply Young's inequality to these four terms above, we obtain 
\begin{align*}
& 2\gamma_k(f(x^{k+1})- f(w^k)+g(x^{k+1})-g(z^k)) \\  
\leq & 2\gamma_k G(2\|x^k-x^{k+1}\|_1+2\|w^k-x^k\|_1+\|z^k-\bar{x}^k\|_1+\|\bar{x}^k-w^k\|_1) \\
\leq & \left(\frac{8\gamma_k^2G^2}{\sigma} + \frac{\sigma}{2}\|x^{k+1}-x^k\|_1^2\right)+\left(\frac{8\gamma_k^2G^2}{\sigma} + \frac{\sigma}{2}\|w^k-x^k\|_1^2\right) \\ & \left(\frac{2\gamma_k^2G^2}{\sigma} + \frac{\sigma}{2}\|z^k-\bar{x}^k\|_1^2\right) + \left(\frac{4\gamma_k^2G^2}{\sigma} + \frac{\sigma}{4}\|\bar{x}^k-w^k\|_1^2\right) \\
\leq & \gamma_k^2 + D_h(x^{k+1},x^k) + D_h(w^k,x^k) + D_h(z^k,\bar{x}^k) + \frac{1}{2}D_h(\bar{x}^k,w^k),
\end{align*}
where in the last inequality we omitted the coefficient of $\gamma_k^2$, and we used the $\sigma$-strong convexity of $h$. Combining this inequality with \eqref{final-bound-inequality}, we get
\[
2\gamma_k(f(x^{k+1})+g(x^{k+1})-f(y)-g(y)) +2D_h(y,x^{k+1}) \leq 2D_h(y,x^k) + \gamma_k^2.
\]
Now summing this over $k=0,1,\ldots,N-1$, we get 
\begin{align*}
\sum_{k=0}^{N-1}\gamma_k (f(x^{k+1})+g(x^{k+1})-f(y)-g(y)) \leq D_h(y,x^0) + \sum_{k=0}^{N-1}\gamma_k^2.
\end{align*}
That is, 
\[
\min_{0\leq k\leq N-1} \left(f(x^{k+1})+g(x^{k+1})\right) - (f(y)+g(y))  \leq \frac{D_h(y,x^0)}{\sum_{k=0}^{N-1}\gamma_k} + \frac{\sum_{k=0}^{N-1}\gamma_k^2}{\sum_{k=0}^{N-1}\gamma_k}.
\]
By choosing $y=x^*$ and using $\gamma_k = \frac{1}{\sqrt{k}}$, we obtain the convergence and the sublinear rate.

\end{document}